\documentclass{amsart}
\usepackage{amsmath,amssymb, amsbsy}
\usepackage[dvips]{graphicx}
\newcommand{\supess}{\mathop{\rm ess\,sup}}
\newcommand{\media}{\mkern12mu\hbox{\vrule height4pt depth-3.2pt
    width6pt} \mkern-17.9mu\int} 
\newcommand{\mediapiccola}{\mkern5mu\hbox{\vrule height4pt depth-3.5pt
    width4pt} \mkern-12.5mu\int}
\newcommand{\R}{{\mathbb R}}
\newcommand{\N}{{\mathbb N}}
\newcommand{\Z}{{\mathbb Z}}

\renewcommand{\a }{\alpha}

\newcommand{\D }{\Delta}
\newcommand{\e }{\varepsilon}

\newcommand{\n }{\nabla}

\newcommand{\Di}{{\mathcal D}^{1,2}(\R^N)}

\renewcommand{\geq }{\geqslant}
\renewcommand{\le }{\leqslant}
\renewcommand{\leq }{\leqslant}

\newenvironment{pf}{\noindent{\sc Proof}.\enspace}{\hfill\qed\medskip}
\newenvironment{pfn}[1]{\noindent{\bf Proof of {#1}.\enspace}}{\hfill\qed\medskip}
\newtheorem{Theorem}{Theorem}[section]

\newtheorem{Lemma}[Theorem]{Lemma}
\newtheorem{Proposition}[Theorem]{Proposition}

\newtheorem{remark}[Theorem]{Remark}
\begin{document}

\title[Schr\"odinger equations with dipole--type 
potentials]{On the behavior of
solutions to Schr\"odinger equations with dipole type potentials near the
singularity}

\author[Veronica Felli]{Veronica Felli}
\address{\hbox{\parbox{5.7in}{\medskip\noindent{Veronica Felli:
        Universit\`a di Milano Bicocca, Dipartimento di Statistica,
        Via Bicocca degli Arcimboldi 8, 20126 Milano, Italy.
        \em{E-mail address: }{\tt veronica.felli@unimib.it.}}}}}
\author[Elsa M. Marchini \and Susanna Terracini]{Elsa M. Marchini \and
  Susanna Terracini}
\address{\hbox{\parbox{5.7in}{\medskip\noindent{Elsa M. Marchini,
        Susanna Terracini: Universit\`a di Milano Bicocca,
        Dipartimento di Ma\-t\-ema\-ti\-ca e Applicazioni, Via Cozzi 53, 20125
        Milano, Italy. \em{E-mail addresses: }{\tt
          elsa.marchini@unimib.it, susanna.terracini@unimib.it}.}}}}

\date{December 21, 2006}

\thanks{First and third author supported by Italy MIUR, national project ``Variational
  Methods and Nonlinear Differential Equations''.
  }

\begin{abstract}
  \noindent Asymptotics of solutions to Schr\"odinger equations with
  singular dipole-type potentials is investigated.  We evaluate the
  exact behavior near the singularity of solutions to elliptic
  equations with potentials which are purely angular multiples of
  radial inverse-square functions. Both the linear and the semilinear
  (critical and subcritical) cases are considered.
\end{abstract}

\maketitle

\begin{center}
{\it Dedicated to Prof. Norman Dancer on the occasion of his 60th birthday.}
\end{center}


\bigskip
\section{Introduction and statement of the main results}\label{intro}

In nonrelativistic molecular physics, the interaction between an
electric charge and the dipole moment ${\mathbf D}\in\R^N$ of a
molecule is described by an inverse square potential with an
anisotropic coupling strength. In particular the Schr\"odinger
equation for the wave function of an electron interacting with a polar
molecule (supposed to be point-like) can be written as
$$
\bigg(\!\!-\frac{\hbar^2}{2m}\,\Delta+e\,\frac{x\cdot {\mathbf
    D}}{|x|^3}-E\bigg)\Psi=0,
$$
where $e$ and $m$ denote respectively the charge and the mass of the
electron and ${\mathbf D}$ is the dipole moment of the molecule, see
\cite{leblond}.

We aim to describe  the asymptotic
behavior near the singularity of solutions to equations 
associated to dipole-type Schr\"odinger operators of the form 
\begin{equation*}
L_{\lambda,{\mathbf d}}:=-\D-\dfrac{\lambda\,(x\cdot{\mathbf d})}{|x|^3}
\end{equation*}
in $\R^N$, where $N\geq 3$, $\lambda=\frac{2me|{\mathbf
    D}|}{\hbar^2}$, being $|{\mathbf D}|$ the magnitude of the dipole
moment ${\mathbf D}$, and ${\mathbf d}={\mathbf D}/{|{\mathbf D}|}$
denotes the orientation of ${\mathbf D}$.  A precise estimate of such
a behavior is the first step in the analysis of fundamental properties
of Schr\"odinger operators, such as positivity, essential
self-adjointness, and spectral features. Such an analysis has been
carried out in  \cite{FMT} for the case of 
Schr\"odinger operators with multipolar Hardy potentials based on the
regularity results proved in \cite{FS3} for degenerate elliptic
equations. In a forthcoming paper, the authors will apply Theorem \ref{t:1} 
to establish spectral properties of multi-singular dipole type operators.

We emphasize that, from the mathematical point of
view, potentials of the form $\frac{\lambda\,(x\cdot{\mathbf d})}{|x|^3}$
have the same order of homogeneity as inverse square potentials and
consequently share many features with them, such as invariance by
scaling and Kelvin transform, as well as no inclusion in the Kato
class.  We mention that Schr\"odinger equations with Hardy-type
singular potentials have been largely studied, see e.g. \cite{AFP,
   egnell,  GP, Jan,  SM, terracini96} and references
therein.

More precisely, in this paper we deal with a more general class of
Schr\"odinger operators including those with dipole-type potentials,
namely with operators whose potentials are purely angular multiples of
radial inverse-square potentials:
\begin{equation*}
{\mathcal L}_{a}:=-\D-\dfrac{a(x/|x|)}{|x|^2}
\end{equation*}
in $\R^N$, where $N\geq 3$ and $a\in L^{\infty}({\mathbb S}^{N-1})$. 

The problem of establishing the asymptotic behavior of solutions to
elliptic equations near an isolated singular point has been studied by
several authors in a variety of contexts, see e.g. \cite{pinchover94}
for Fuchsian type elliptic operators and \cite{DV2} for Fuchsian type
weighted operators. The asymptotics we derive in this work is not
contained in the aforementioned papers, which prove the existence of
the limit at the singularity of any quotient of two positive
solutions in some linear and semilinear cases which, however, do not not include 
the perturbed linear case and the critical nonlinear case treated here. Moreover,
besides proving the existence of such a limit, we also obtain a Cauchy type representation 
formula for it, see (\ref{eq:54}) and (\ref{eq:4}). 

We also quote \cite{murata}, where
asymptotics at infinity is established for perturbed
inverse square potentials and in  some particular  nonradial case.
 H\"older continuity results for degenerate elliptic equations with
singular weights and  asymptotic analysis of the behavior of solutions 
near the pole  are contained in \cite{FS3}.

As a natural setting to study the properties of operators
${\mathcal L}_{a}$, we introduce the functional space $\Di$ defined
as the completion of $C^\infty_{\rm c}(\R^N)$ with respect to the
Dirichlet norm
$$
\|u\|_{\Di}:=\bigg(\int_{\R^N}|\n u(x)|^2\,dx\bigg)^{\!\!1/2}.
$$ 
In order to discuss the positivity properties of the Schr\"odinger
operator ${\mathcal L}_{a}$ in $\Di$, we consider the best constant in
the associated Hardy-type inequality
\begin{equation}
\label{eq:bound}
\Lambda_N(a):=\sup_{u\in\Di\setminus\{0\}}\dfrac{{\displaystyle
{\int_{\R^N}{{|x|^{-2}}{a(x/|x|)}\,u^2(x)\,dx}}}}
{{\displaystyle{\int_{\R^N}{|\n u(x)|^2\,dx}}}}.
\end{equation} 
By the classical Hardy inequality,  $\Lambda_N(a)\leq
\frac4{(N-2)^2}\supess_{\mathbb S^{N-1}}a$, where $\supess_{\mathbb
  S^{N-1}}a$ denotes the essential supremum of $a$ in $\mathbb
S^{N-1}$.  We also notice that, in the dipole case, by rotation invariance,
$\Lambda_N\big(\frac{\lambda\,x}{|x|}\cdot{\mathbf d}\big)$ does not depend on
${\mathbf d}$.

It is easy to verify that the quadratic form
associated to ${\mathcal L}_{a}$ is positive definite in $\Di$ if and
only if $\Lambda_N(a)<1$.  The relation between the value
$\Lambda_N(a)$ and the first eigenvalue of the angular component of
the operator on the sphere $\mathbb S^{N-1}$ is discussed in section
\ref{sec:spectr-angul-comp}, see Lemmas \ref{l:ln} and \ref{l:pos}.
More precisely, Lemma \ref{l:pos} ensures that the quadratic form
associated to ${\mathcal L}_{a}$ is positive definite if and only if
$$
\mu_1>-\bigg(\frac{N-2}2\bigg)^{\!\!2},
$$
where $\mu_1=\mu_1(a,N)$ is the first eigenvalue of the
operator $-\D_{\mathbb S^{N-1}}-a(\theta)$ on $\mathbb S^{N-1}$ (see
Lemma \ref{l:spe}). We  denote  by $\psi_1$ the associated positive
$L^2$-normalized eigenfunction and  set 
\begin{equation}\label{eq:sigma}
\sigma=\sigma(a,N):=-\frac{N-2}2+\sqrt{\bigg(\frac{N-2}2\bigg)^{\!\!2}+\mu_1}.
\end{equation}
In the spirit of the well-known Riemann removable singularity theorem, we describe
 the behavior of solutions to linear
Schr\"odinger equations with a dipole-type singularity localized in a
neighborhood of $0$.
\begin{Theorem}\label{t:1}
  Let $\Omega\subset\R^N$ be a bounded open set such that $0\in\Omega$
  and let $q\in L^{\infty}_{\rm loc}(\Omega\setminus\{0\})$ be such that
$q(x)=O(|x|^{-(2-\e)})$ as $|x|\to0$ for some $\e>0$. Assume that
  $a\in L^{\infty}({\mathbb S}^{N-1})$ satisfies $\Lambda_N(a)<1$ and
  $u\in H^1(\Omega)$, $u\geq 0$ a.e. in $\Omega$, $u\not\equiv 0$,
  weakly solves
\begin{equation}\label{eq:45}
  -\Delta u(x)-\frac{a(x/|x|)}{|x|^2}\,u(x)
=q(x)\,u(x)\quad\text{ in }\Omega,
\end{equation}
i.e.
$$
\int_{\Omega}\n u(x)\cdot\n w(x)\,dx-\int_{\Omega} \frac{a(x/|x|)}{|x|^2}
\,u(x)w(x)\,dx=\int_{\Omega}
q(x)u(x)w(x)\,dx,\quad\text{for all }w\in
H^1_0(\Omega).
$$
 Then the function 
$$
x\mapsto \frac{u(x)}{|x|^{\sigma}\psi_1(x/|x|)}
$$
is continuous in $\Omega$ and 
\begin{align}\label{eq:54}
\lim_{|x|\to 0}\frac{u(x)}{|x|^{\sigma}\psi_1\big(\frac x{|x|}\big)}=
\int_{{\mathbb
S}^{N-1}}\bigg(R^{-\sigma}&u(R\theta)+\int_0^R 
{{\frac{s^{1-\sigma}}
{2\sigma+N-2}}}\,q(s\,\theta)u(s\,\theta)\,ds
\\
&\notag-R^{-2\sigma-N+2}\int_0^R 
{{\frac{s^{N-1+\sigma}}
  {2\sigma+N-2}}}\,q(s\,\theta)u(s\,\theta)\,ds
\bigg)\psi_1(\theta)
  \,dV(\theta),
\end{align}
for all $R>0$ such that $\overline{B(0,R)}:=\{x\in\R^N:\ |x|\leq
R\}\subset\Omega$.
\end{Theorem}

We notice that \eqref{eq:54} is actually a {\it Cauchy's integral type formula}
for $u$. Moreover the term at the right hand side is independent of $R$.
 In the case in which the perturbation $q$ is radial
then an analogous formula holds also for changing sign solutions to \eqref{eq:45}, 
see Remark \ref{r:changingsign}.

\smallskip
If the perturbing potential $q$ satisfies some proper summability
condition, instead of the stronger control on the blow-up rate at the
singularity required in Theorem \ref{t:1}, a Brezis-Kato type
argument, see \cite{BrezisKato}, allows us to derive an
upper estimate on the behavior of solutions. For any $q\geq 1$, we denote as
$L^q(\varphi^{2^*}\!\!,\Omega)$ the weighted $L^q$-space endowed with
the norm
$$
\|u\|_{L^q(\varphi^{2^*}\!\!,\Omega)}:=
\bigg(\int_{\Omega}\varphi^{2^*}(x)|u(x)|^s\,dx\bigg)^{\!\!1/q},
$$
where $2^*=\frac{2N}{N-2}$ is the critical Sobolev exponent and
 and $\varphi$ denotes the weight function
\begin{equation}\label{eq:weight}
\varphi(x):=|x|^{\sigma}\psi_1(x/|x|). 
\end{equation}
The following
Brezis-Kato type result holds.
\begin{Theorem}\label{t:bk}
  Let $\Omega\subset\R^N$ be a bounded open set such that $0\in\Omega$,
    $a\in L^{\infty}({\mathbb S}^{N-1})$ satisfying $\Lambda_N(a)<1$, and $V\in
  L^s(\varphi^{2^*}\!\!,\Omega)$ for some $s>N/2$. Then, for any
  $\Omega' \Subset \Omega$, there exists a positive constant
  $$
C=C\big(N,a,\|V\|_{L^s(\varphi^{2^*}\!\!,\Omega)},\mathop{\rm
    dist} (\Omega',\partial\Omega),\mathop{\rm diam}\Omega \big)
$$
  depending only on $N$, $a$,
  $\|V\|_{L^s(\varphi^{2^*}\!\!,\Omega)}$, $\mathop{\rm dist}
  (\Omega',\partial\Omega)$, and $\mathop{\rm diam}\Omega$, such that
for any weak
  $H^1(\Omega)$-solution  $u$ of
\begin{equation}\label{eq:67}
  -\Delta u(x)-\frac{a(x/|x|)}{|x|^2}\,u(x)=\varphi^{2^*-2}(x)V(x)u(x),
\quad \text{in }\Omega,
\end{equation}
i.e. for any $u\in H^1(\Omega)$ satisfying
$$
\int_{\Omega}\n u(x)\cdot\n w(x)\,dx-\int_{\Omega}
\frac{a(x/|x|)}{|x|^2}\,u(x)w(x)\,dx=\int_{\Omega}
\varphi^{2^*-2}(x)V(x)u(x)w(x)\,dx,
$$
for all $w\in H^1_0(\Omega)$, there holds $\frac {u}{\varphi}\in
L^{\infty}(\Omega')$ and
\begin{equation*}
  \bigg\|\frac{u}{\varphi}\bigg\|_{L^{\infty}(\Omega')}\leq 
C\,\|u\|_{L^{2^*}(\Omega)}.
\end{equation*}
\end{Theorem}

The Brezis-Kato procedure can be applied also to semilinear problems
with at most critical growth, thus providing an upper bound for
solutions and then reducing the semilinear problem to a linear one
with enough control on the potential at the singularity to apply
Theorem \ref{t:1} and to recover the exact asymptotic behavior.

\begin{Theorem}\label{t:semi}
  Let $\Omega$ be a bounded open set containing $0$,
$a\in L^{\infty}({\mathbb S}^{N-1})$ satisfying $\Lambda_N(a)<1$, 
and $f:\Omega\times\R\to\R$ such that,
for some positive constant $C$, 
$$
\bigg|\frac{f(x,u)}u\bigg|\leq C\,\big(1+|u|^{2^*-2}\big)\quad\text{for a.e. }
(x,u)\in \Omega\times\R.
$$ 
Assume that
  $u\in H^1(\Omega)$, $u\geq 0$ a.e. in $\Omega$, $u\not\equiv 0$,
  weakly solves 
\begin{equation}\label{eq:72}
-\Delta u(x)-\frac{a(x/|x|)}{|x|^2}\,u(x)=f(x,u(x)),\quad \text{in }\Omega,
\end{equation}
i.e. 
$$
\int_{\Omega}\n u(x)\cdot\n w(x)\,dx-\int_{\Omega} \frac{a(x/|x|)}{|x|^2}
\,u(x)w(x)\,dx=\int_{\Omega}
f(x,u(x))w(x)\,dx,\quad\text{for all }w\in
H^1_0(\Omega).
$$
 Then the function 
$$
x\mapsto \frac{u(x)}{|x|^{\sigma}\psi_1(x/|x|)}
$$
is continuous in $\Omega$ and 
\begin{align}\label{eq:4}
\lim_{|x|\to 0}\frac{u(x)}{|x|^{\sigma}\psi_1\big(\frac x{|x|}\big)}=
\int_{{\mathbb
S}^{N-1}}\bigg(r^{-\sigma}&u(r\theta)+\int_0^r 
{{\frac{s^{1-\sigma}}
{2\sigma+N-2}}}\,f\big(s\,\theta,u(s\,\theta)\big)\,ds
\\
&\notag-r^{-2\sigma-N+2}\int_0^r 
{{\frac{s^{N-1+\sigma}}
  {2\sigma+N-2}}}\,f\big(s\,\theta,u(s\,\theta)\big)\,ds
\bigg)\psi_1(\theta)
  \,dV(\theta),
\end{align}
for all $r>0$ such that $\overline{B(0,r)}:=\{x\in\R^N:\ |x|\leq r\}\subset\Omega$.
\end{Theorem}

\goodbreak

\medskip
\noindent
{\bf Notation. } We list below some notation used throughout the
paper.\par
\begin{itemize}
\item[-]$B(a,r)$ denotes the ball $\{x\in\R^N: |x-a|<r\}$ in $\R^N$ with
center at $a$ and radius $r$.
\item[-] $dV$ denotes the volume element on the sphere ${\mathbb S}^{N-1}$.
\item[-] $\omega_N$ denotes the volume of the unit sphere ${\mathbb
    S}^{N-1}$, i.e. $\omega_N=\int_{{\mathbb S}^{N-1}}dV(\theta)$.
\item[-] for any $a\in L^1({\mathbb S}^{N-1})$, $\mediapiccola_{{\mathbb
    S}^{N-1}}
  a(\theta)\,dV(\theta)$ denotes the mean of $a$ on ${\mathbb
    S}^{N-1}$, i.e. 
$$\media_{{\mathbb
    S}^{N-1}}
  a(\theta)\,dV(\theta)=\frac1{\omega_N}\int_{{\mathbb
    S}^{N-1}} a(\theta)\,dV(\theta).
$$
\item[-] the symbol $\supess$ stands for essential supremum. 
\end{itemize}

\section{Spectrum of the angular component}\label{sec:spectr-angul-comp}

Due to the structure of the dipole-type potential of equation
(\ref{eq:45}), a natural approach to describe the solutions seems to
be the separation of variables. To employ such a technique, we need,
as a starting point, the description of the spectrum of the angular
part of dipole Schr\"odinger operators.

\begin{Lemma}\label{l:spe}
Let $a\in L^{\infty}({\mathbb S}^{N-1})$. Then the operator $-\D_{\mathbb
S^{N-1}}-a(\theta)$ on $\mathbb S^{N-1}$ admits
a diverging sequence of eigenvalues
$\mu_1<\mu_2\leq\cdots\leq\mu_k<\cdots$ the first of which has the
following properties:
\begin{itemize}
\item[(i)] $\mu_1$ is  simple;\\[-7pt]
\item[(ii)] $\mu_1$ can be characterized as
\begin{equation}\label{firsteig}
\mu_1=\min_{\psi\in H^1(\mathbb
S^{N-1})\setminus\{0\}}\frac{\int_{\mathbb S^{N-1}}|\n_{\mathbb
S^{N-1}}\psi(\theta)|^2\,dV(\theta)-\int_{\mathbb S^{N-1}}a(\theta)
\psi^2(\theta)\,dV(\theta)}{\int_{\mathbb S^{N-1}}\psi^2(\theta)\,dV(\theta)};
\end{equation}
\item[(iii)] $\mu_1$ is attained by  a $C^1$  positive eigenfunction $\psi_1$
such that $\min_{\mathbb S^{N-1}}\psi_1>0$;\\[-7pt]
\item[(iv)] if, for some $\kappa\in\R$, $a(\theta)=\kappa$ for a.e.
  $\theta\in \mathbb S^{N-1}$, then  $\mu_1=-\kappa$.\\[-7pt]
\item[(v)] if $a$ is not constant, then $-\supess_{\mathbb S^{N-1}}
  a<\mu_1<-\mediapiccola_{{\mathbb S}^{N-1}} a(\theta)\,dV(\theta)$.
\end{itemize}
\end{Lemma}

\begin{pf}
We prove assertion (v), being (i), (ii), (iii), and (iv) quite standard.

Since the function $\psi\equiv1$ satisfies 
$$
\int_{\mathbb S^{N-1}}|\n_{\mathbb
S^{N-1}}\psi(\theta)|^2\,dV(\theta)-\int_{\mathbb
S^{N-1}}a(\theta)\psi^2(\theta)\,dV(\theta)=-\int_{\mathbb
S^{N-1}}a(\theta)\,dV(\theta),
$$ 
we deduce that $\mu_1\leq -\mediapiccola_{{\mathbb S}^{N-1}}
a(\theta)\,dV(\theta)$.  In order to prove the strict inequality, we
argue by contradiction and assume that $\mu_1=-\mediapiccola_{{\mathbb
    S}^{N-1}} a(\theta)\,dV(\theta)$. Then $\psi_1\equiv1$ attains the
minimum value in \eqref{firsteig} but, since $a$ is not constant, it
does not satisfy equation $-\D_{\mathbb
  S^{N-1}}\psi_1-a(\theta)\psi_1=0$, a contradiction. We can thereby
conclude that $\mu_1< -\mediapiccola_{{\mathbb S}^{N-1}}
a(\theta)\,dV(\theta)$.

From \eqref{firsteig}, \cite[Lemma 1.1]{terracini96}, and the optimality of the
best constant in Hardy's inequality, it follows that
\begin{align*}
\mu_1&>\inf_{\psi\in H^1(\mathbb
S^{N-1})\setminus\{0\}}\frac{\int_{\mathbb S^{N-1}}|\n_{\mathbb
S^{N-1}}\psi(\theta)|^2\,dV(\theta)-\supess_{\mathbb S^{N-1}}
  a\int_{\mathbb S^{N-1}}
\psi^2(\theta)\,dV(\theta)}{\int_{\mathbb S^{N-1}}\psi^2(\theta)\,dV(\theta)}\\[5pt]
&=\inf_{u\in\Di\setminus\{0\}}\frac{\int_{\mathbb R^N}|\n
u(x)|^2\,dx-\Big(\big(\frac{N-2}{2}\big)^2+\supess_{\mathbb S^{N-1}}
  a\Big)\int_{\mathbb
R^N}\frac{u^2(x)}{|x|^2}\,dx}{\int_{\mathbb
R^N}\frac{u^2(x)}{|x|^2}\,dx}=-\supess_{\mathbb S^{N-1}}
  a,
\end{align*} 
thus proving the left part of the inequality
stated in (v).
\end{pf}

\noindent 
The asymptotic behavior of eigenvalues $\mu_k$ as $k\to+\infty$ is
described by Weyl's law, which is recalled in the theorem below. We
refer to \cite{reedsimon4, SafarovVassiliev} for a proof.
\begin{Theorem}[\bf Weyl's law]\label{t:weyl}
For  $a\in L^{\infty}({\mathbb S}^{N-1})$, let $\{\mu_k\}_{k\geq1}$ be the eigenvalues of
the operator $-\D_{\mathbb S^{N-1}}-a(\theta)$ on
$\mathbb S^{N-1}$. Then
\begin{equation}\label{eq:weyl}
\mu_k=C(N,a)k^{2/(N-1)}\big(1+o(1)\big)\quad\text{as }k\to+\infty,
\end{equation}
for some positive constant $C(N,a)$ depending only on $N$ and $a$.
\end{Theorem}

The following lemma provides an estimate of the $L^{\infty}$-norm of
eigenfunctions of the operator $-\D_{\mathbb
S^{N-1}}-a(\theta)$ in terms of the corresponding eigenvalues.
For classical results about $L^{\infty}$ estimates of eigenfunctions of 
Schr\"odinger operators we refer to \cite[\S 4]{simon2000} and references therein.
\begin{Lemma}\label{l:stiautf}
For  $a\in L^{\infty}({\mathbb S}^{N-1})$ and $k\in\N\setminus\{0\}$, 
let $\psi_k$ be a $L^2$-normalized
eigenfunction of the Schr\"odinger operator $-\D_{\mathbb
S^{N-1}}-a(\theta)$ on the sphere associated to
the $k$-th eigenvalue $\mu_k$, i.e.
\begin{equation}
\label{eq:2rad}
\begin{cases}
-\D_{{\mathbb S}^{N-1}}\psi_k(\theta)-a(\theta)\,\psi_k(\theta)
=\mu_k\,\psi_k(\theta),&\text{in }{\mathbb S}^{N-1},\\[3pt]
\int_{{\mathbb S}^{N-1}}|\psi_k(\theta)|^2\,dV(\theta)=1.
\end{cases}
\end{equation}
Then, there exists a constant $C_1$ depending only on $N$ 
and $a$ such that 
$$
|\psi_k(\theta)|\leq C_1\, |\mu_k|^{\lfloor (N-1)/4\rfloor+1},
$$
where $\lfloor \cdot\rfloor$ denotes the floor function, i.e. $\lfloor
x\rfloor:=\min\{j\in\Z:\ j\leq x\}$.
\end{Lemma}

\begin{pf}
Using classical elliptic regularity theory and bootstrap methods, we
can easily prove that for any $j\in\N$ there exists a constant
$C(N,j)$, depending only on $j$ and $N$, such that
$$
\|\psi_k\|_{W^{2,\frac{2(N-1)}{(N-1)-4(j-1)}}({\mathbb S}^{N-1})}\leq
C(N,j)\big(\mu_k+\|a\|_{L^{\infty}({\mathbb S}^{N-1})}\big)^j.
$$
Choosing $j=\big\lfloor \frac{N-1}4\big\rfloor+1$, by Sobolev's inclusions we
deduce that
$$
W^{2,\frac{2(N-1)}{(N-1)-4(j-1)}}({\mathbb
S}^{N-1})\hookrightarrow C^{0,\alpha}({\mathbb
S}^{N-1})\hookrightarrow L^{\infty}({\mathbb S}^{N-1}),
$$
for any $0<\alpha<2\big(1-\frac {N-1}4+\big\lfloor\frac
{N-1}4\big\rfloor\big)$, thus implying the required estimate.
\end{pf}

\noindent

Arguing as in the proof of \cite[Lemma
1.1]{terracini96}, we can deduce the following characterization of
$\Lambda_N(a)$.
\begin{Lemma}\label{l:ln}
  For $a\in L^{\infty}({\mathbb S}^{N-1})$, let $\Lambda_N(a)$ be defined
  in (\ref{eq:bound}). Then
\begin{equation}\label{eq:57}
\Lambda_N(a):=\max_{\psi\in H^1(\mathbb
  S^{N-1})\setminus\{0\}}\frac{\int_{\mathbb
    S^{N-1}}a(\theta)\,\psi^2(\theta)\,dV(\theta)}{\int_{\mathbb
    S^{N-1}}|\n_{\mathbb
    S^{N-1}}\psi(\theta)|^2\,dV(\theta)+\big(\frac{N-2}{2}\big)^2\int_{\mathbb
    S^{N-1}} \psi^2(\theta)\,dV(\theta)}.
\end{equation}
\end{Lemma}

We notice that the supremum in (\ref{eq:57}) is achieved due to the
compactness of the embedding $H^1(\mathbb S^{N-1})\hookrightarrow
L^2(\mathbb S^{N-1})$. As a direct consequence of the above lemma, it is possible 
to compare $\Lambda_N(a)$ with the best constant in Hardy's inequality
$$
\frac4{(N-2)^2}=\sup_{u\in\Di\setminus\{0\}}\frac{{
{\int_{\R^N}{{u^2(x)}{|x|^{-2}}\,dx}}}}{{{\int_{\R^N}{|\n u(x)|^2\,dx}}}}.
$$
 Indeed, if $a$ is not constant, there holds
$$
\frac{4}{(N-2)^2}\media_{{\mathbb S}^{N-1}}
a(\theta)\,dV(\theta)
<\Lambda_N(a)<\supess_{{\mathbb S}^{N-1}}a \, \frac4{(N-2)^2},
$$
whereas, if $a(\theta)=\kappa$ for a.e.
  $\theta\in \mathbb S^{N-1}$ and for some $\kappa\in\R$, then
$\Lambda_N(a)={4\kappa}/{(N-2)^2}$.

Let us consider the quadratic form associated to the Schr\"odinger
operator ${\mathcal L}_{a}$, i.e. 
$$
Q_{a}(u):=\int_{\R^N} |\n  u(x)|^2dx-
\int_{\R^N}\frac{a(x/|x|)\,u^2(x)}{|x|^2}\,dx.
$$
The problem of positivity of $Q_{a}$ is solved in the following lemma.
\begin{Lemma}\label{l:pos}
Let $a\in L^{\infty}({\mathbb S}^{N-1})$. The following conditions are equivalent:
\begin{align*} {\rm i)}\quad&Q_{a} \text{ is positive definite, i.e. }
  \inf_{u\in\Di\setminus\{0\}}\frac{Q_{a}(u)}
  {\int_{\R^N}|\n u(x)|^2\,dx}>0; \\
  {\rm ii)} \quad&\Lambda_N(a)<1;\\
  {\rm iii)}\quad&
  \mu_1>-\big({\textstyle{\frac{N-2}2}}\big)^{\!2}\text{ where $\mu_1$
    is defined in (\ref{firsteig})}.
\end{align*}
\end{Lemma}
\begin{pf}
  The equivalence between i) and ii) follows from the definition of
  $\Lambda_N(a)$, see (\ref{eq:bound}). On the other hand, \cite[Proposition
  1.3 and Lemma 1.1]{terracini96} ensure that i) is equivalent to iii).
\end{pf}

\begin{remark}\label{r:posl}
  We notice that in the case of dipole potentials, namely if
  $a(\theta)=\lambda \, \frac{x\cdot {\mathbf d}}{|x|}$, then, rewriting
  (\ref{eq:57}) in spherical coordinates and exploiting the symmetry
  with respect to the dipole axis, $\Lambda_N(a)$ can be characterized
  as
$$
\Lambda_N\Big(\lambda \, \frac{x\cdot {\mathbf d}}{|x|}\Big)=\lambda \, \sup_{w\in H^1_0(0,\pi)}
\frac{\int_0^\pi \cos\theta\,w^2(\theta)\,d\theta}
{\int_0^\pi\Big[|w'(\theta)|^2+\frac{(N-2)(N-4)}4
(\sin\theta)^{-2}w^2(\theta)\Big]\,d\theta}.
$$
In dimension $N=3$, a Taylor's expansion of 
$\Lambda_N\big(\lambda \, \frac{x\cdot {\mathbf d}}{|x|}\big)$ 
near $\lambda=0$ can be found in \cite{leblond}. 

There is no explicit formula for the values $\Lambda_N\big(\lambda \,
\frac{x\cdot {\mathbf d}}{|x|}\big)$. When $N=4$, it can be expressed in terms
of Mathieu's special functions.
The results of a numerical approximation of $\Lambda_N\big(\lambda \,
\frac{x\cdot {\mathbf d}}{|x|}\big)$  performed with
a finite difference method are listed in  table \ref{tab:hardyconstants}, which
highlights how the dipole Hardy-type constant detaches the classical
Hardy constant more and more as the dimension grows.

\small\rm
\bigskip\noindent
\begin{table}[h]
\begin{tabular}{|l|l|l|}
  \hline
   & & \\[-5pt]
 $N$ &$\big(\Lambda_N(1)\big)^{-1}=\frac{(N-2)^2}4$ & $\Big
[\Lambda_N\big( \frac{x\cdot {\mathbf d}}{|x|}\big)\Big]^{-1}$ \\[7pt] \hline
  & & \\[-5pt]
  3& 0.25& 1.6398  \\[3pt]
  \hline
  & & \\[-5pt]
  4 & 1&3.7891
  \\[3pt]
  \hline
  & & \\[-5pt]
  5& 2.25& 7.5831   \\[3pt]
  \hline
  & & \\[-5pt]
  6& 4&  12.6713 \\[3pt]
  \hline
    \end{tabular}\quad
\begin{tabular}{|l|l|l|}
\hline
   & & \\[-5pt]
 $N$ &$\big(\Lambda_N(1)\big)^{-1}=\frac{(N-2)^2}4$ & $\Big[\Lambda_N
\big( \frac{x\cdot {\mathbf d}}{|x|}\big)\Big]^{-1}$ \\[7pt] \hline
& & \\[-5pt]
  7& 6.25&  19.0569 \\[3pt]
  \hline
  & & \\[-5pt]
  8& 9 &  26.7407 \\[3pt]
  \hline
  & & \\[-5pt]
  9& 12.25 &  35.7231  \\[3pt]
  \hline
  & & \\[-5pt]
  10& 16&  46.0044 \\[3pt]
  \hline
   \end{tabular}
\\[8pt]
\caption{\scriptsize Some numerical approximations of 
  $\Lambda_N \big( \frac{x\cdot {\mathbf d}}{|x|}\big)$  
obtained by finite difference
  with 10000 steps.}
\label{tab:hardyconstants}\end{table}
\normalsize
\end{remark}

\section{A Brezis-Kato type lemma}\label{sec:bk}

In this section, we follow the procedure developed by Brezis and Kato
in \cite{BrezisKato} to control from above the behavior of solutions
to Schr\"odinger equations with dipole type potentials, in order to
prove Theorem \ref{t:bk}.

 Let us consider the weight $\varphi$ introduced in (\ref{eq:weight}) and define 
the  weighted $H^1$-space $H^1_{\varphi}(\Omega)$
 as the closure of $C^\infty(\bar \Omega)$ with respect to
\begin{align*}
  \|u\|^2_{H^1_{\varphi}(\Omega)} := \int_\Omega \varphi^2(x)
  \left(|\nabla u|^2+|u|^2\right)\,dx, 
\end{align*}
and the space ${\mathcal D}^{1,2}_{\varphi}(\Omega)$
as the closure of $C^\infty_{\rm c}(\Omega)$ with respect to
\begin{align*}
  \|u\|^2_{{\mathcal D}^{1,2}_{\varphi}(\Omega)} := \int_\Omega \varphi^2(x)
  |\nabla u|^2\,dx.
\end{align*}
By the Caffarelli-Kohn-Nirenberg inequality (see \cite{CKN} and
\cite{CatrinaWang}) and the definition of $\varphi$, it follows that,
for any $w\in {\mathcal D}^{1,2}_{\varphi}(\Omega)$,
\begin{align}
\label{eq:63}
\left(\int_{\Omega} \varphi^{2^*}(x)|w(x)|^{2^*}\,dx\right)^{2/{2^*}}
\leq {\mathcal C}_{N,a}\int_{\Omega}
  \varphi^{2}(x)|\nabla w|^2 \,dx,
\end{align}
for some positive constant ${\mathcal C}_{N,a}$ depending only
on $N$ and $a$.

\begin{Lemma}\label{l:bk1}
    Let $\Omega$ be a bounded open set containing $0$ and
 $v\in H^1_{\varphi}(\Omega)\cap L^q(\varphi^{2^*}\!\!,\Omega)$, $q>1$, be a
  weak solution to
\begin{equation}\label{eq:65}
-\mathop{\rm div}(\varphi^2(x)\n
v(x))=\varphi^{2^*}(x)V(x)v(x),\quad\text{in }\Omega,
\end{equation}
where $V\in L^s(\varphi^{2^*}\!\!,\Omega)$ for some $s>\frac N2$. Then,
there exists a positive constant $\widetilde C=\widetilde
C(a,N)$ depending only on $a$ and  $N$,  such that for
any $\Omega' \Subset \Omega$, $v\in
L^{\frac{2^*q}2}(\varphi^{2^*}\!\!,\Omega') $ and
\begin{multline}\label{eq:66}
  \|v\|_{L^{\frac{2^*q}2}(\varphi^{2^*}\!\!,\Omega')}\leq 
\widetilde C^{\frac1q} 
\big(\!\mathop{\rm diam}\Omega\big)^{\frac{\sigma(2-2^*)}q}\times\\
\quad \times\bigg(\frac{8}{C(q)}\frac{4}{(\mathop{\rm dist}
    (\Omega',\partial\Omega))^2}+
  \frac{4(q+2)}{(\mathop{\rm dist}(\Omega',\partial\Omega))^2}
  +\frac{2\ell_q}{C(q)}\bigg)^{\!\!\frac1q} \|v\|_{L^{q}(\varphi^{2^*}\!\!,\Omega)},
\end{multline}
where $C(q):=\min\big\{\frac14, \frac4{q+4}\big\}$ and 
$$
\ell_q=\bigg[\max\bigg\{8\,{\mathcal
  C}_{N,a}\|V\|_{L^s(\varphi^{2^*}\!\!,\Omega)}^{2s/N},
\frac{q+4}2{\mathcal
  C}_{N,a}\|V\|_{L^s(\varphi^{2^*}\!\!,\Omega)}^{2s/N}\bigg\}\bigg]^{\frac{N}{2s-N}}.
$$
\end{Lemma}

\begin{pf}
  H\"older's inequality and \eqref{eq:63} yield for any $w\in{\mathcal
    D}^{1,2}_{\varphi}(\Omega)$
\begin{align}
\label{eq:4bk} \int_{\Omega}&\varphi^{2^*}(x)|V(x)|w^2(x)\,dx
  \leq\ell_q\!\!\!\!\!\!\int\limits_{|V(x)|\leq\ell_q}\!\!\!\!\!\!\varphi^{2^*}(x) w^2(x)\,dx
  +\!\!\!\!\!\!
\int\limits_{|V(x)|\geq \ell_q}\!\!\!\!\!\!\varphi^{\frac4{N-2}}(x)|V(x)|\varphi^{2}(x)w^2(x)\,dx\\
  \notag &\leq \ell_q\int_{\Omega}\varphi^{2^*}(x)w^2(x)\,dx
  +\bigg(\int_{\Omega}\varphi^{2^*}(x)w^{2^*}(x)\,dx\bigg)^{\frac2{2^*}} \bigg(
\int\limits_{|V(x)|\geq \ell_q} \varphi^{2^*}(x)|V(x)|^{\frac N2}\,dx\bigg)^{\frac2N}\\
\notag
&\leq \ell_q\int_{\Omega}\varphi^{2^*}(x)w^2(x)\,dx+C_{N,a} \bigg(\int_{\Omega}
  \varphi^{2}(x)|\nabla w|^2 \,dx\bigg) \bigg(
\int\limits_{|V(x)|\geq \ell_q} \varphi^{2^*}(x)|V(x)|^{\frac N2}\,dx\bigg)^{\frac2N}.
\end{align}
By H\"older's inequality and by the choice of
$\ell_q$ it follows that
\begin{align}\label{eq:70}
  \int\limits_{|V(x)|\geq \ell_q} \varphi^{2^*}(x)&|V(x)|^{\frac
    N2}\,dx \leq
  \bigg(\int_{\Omega}\varphi^{2^*}(x)|V(x)|^s\,dx\bigg)^{\frac{N}{2s}}\bigg(
  \int\limits_{|V(x)|\geq \ell_q} \varphi^{2^*}(x)\,dx\bigg)^{\frac{2s-N}{2s}}\\
  \notag&\leq \bigg(\int_{\Omega}\varphi^{2^*}(x)|V(x)|^s\,dx\bigg)^{\frac{N}{2s}}\bigg(
  \int\limits_{|V(x)|\geq \ell_q}\bigg(\frac{|V(x)|}{\ell_q} \bigg)^s
 \varphi^{2^*}(x)\,dx\bigg)^{\frac{2s-N}{2s}}\\
\notag&\leq \|V\|_{L^s(\varphi^{2^*}\!\!,\Omega)}^{s}\,\ell_q^{-s+\frac
    N2}\leq\min\bigg\{\frac{C_{N,a}^{-1}}{8},
  \frac{2C_{N,a}^{-1}}{q+4}\bigg\}^{\frac N2}, 
\end{align}
and hence from \eqref{eq:4bk} we obtain that
for any $w\in{\mathcal
    D}^{1,2}_{\varphi}(\Omega)$
\begin{align}\label{eq:64}
\int_{\Omega}\varphi^{2^*}(x)|V(x)|w^2(x)\,dx
  \leq\ell_q\int_{\Omega}\varphi^{2^*}(x)w^2(x)\,dx+
\min\bigg\{\frac{1}{8},
  \frac{2}{q+4}\bigg\}\bigg(\int_{\Omega}
  \varphi^{2}(x)|\nabla w|^2 \,dx\bigg).
\end{align}
Let $\eta$ be a nonnegative cut-off function such that 
$$
\mathop{\rm supp}(\eta)
\Subset \Omega, \quad \eta \equiv 1 \text{ on } 
\Omega',\text{ and } |\nabla \eta(x)|\leq \frac{2}{\mathop{\rm dist}
    (\Omega',\partial\Omega)}.
$$
  Set $v^n :=
\min(n,|v|) \in {\mathcal D}^{1,2}_{\varphi}(\Omega)$ and test (\ref{eq:65}) with
$v(v^n)^{q-2}\eta^2 \in {\mathcal D}^{1,2}_{\varphi}(\Omega)$. This leads to
\begin{align*}
(q&-2)\int_{\Omega}\varphi^{2}(x)\eta^2(x)|\n v^n(x)|^2(v^n(x))^{q-2}\,dx+ 
\int_{\Omega}\varphi^{2}(x)\eta^2(x) (v^n(x))^{q-2}|\n v(x)|^2\,dx\\   
&=\int_{\Omega}\varphi^{2^*}(x)V(x)\eta^2(x) v^2(x)(v^n(x))^{q-2}\,dx
-2\int_{\Omega}\varphi^{2}(x) \eta(x) v(x)(v^n(x))^{q-2}\n v(x)\cdot    \n \eta(x)\,dx .
\end{align*}
We use the elementary inequality $2ab \le 1/2 a^2 + 4 b^2$ and obtain
\begin{align}
\label{eq:19bk}
(q&-2)\int_{\Omega}\varphi^{2}(x)\eta^2(x)|\n v^n(x)|^2(v^n(x))^{q-2}\,dx+ 
\frac 12\int_{\Omega}\varphi^{2}(x)\eta^2(x) (v^n(x))^{q-2}|\n v(x)|^2\,dx\\
\notag   
&\leq \int_{\Omega}\varphi^{2^*}(x)V(x)\eta^2(x) v^2(x)(v^n(x))^{q-2}\,dx
+ 4 \int_{\Omega} \varphi^{2}(x)|\n \eta(x)|^2 v^2(x)(v_n(x))^{q-2}\,dx.
\end{align}
Furthermore, an explicit calculation gives
\begin{align}\label{eq:21bk}
\begin{split}
\big|\n \big((v^n)^{\frac q2-1}v\eta)\big|^2
&\le \frac{(q+4)(q-2)}4(v^n)^{q-2}\eta^2|\n v^n|^2+2(v^n)^{q-2}
|\n v|^2\eta^2\\
&\quad+2 (v^n)^{q-2} v^2|\n \eta|^2+\frac{q-2}2 (v^n)^{q}|\n \eta|^2.
\end{split}
\end{align}
Letting $C(q):=\min\big\{\frac14, \frac4{q+4}\big\}$, from
(\ref{eq:19bk}) and (\ref{eq:21bk}) we get
\begin{align}\label{eq:22bk}
  C(q)&\int_{\Omega}\varphi^{2}(x) \big|\n \big((v^n)^{\frac q2-1}
  v\eta\big)(x)\big| ^2\,dx \leq 2(2+C(q))
  \int_{\Omega}\varphi^{2}(x)(v^n(x))^{q-2}v(x)^2|\n \eta(x)|^2 \,dx\\
  \notag &+ C(q)\frac{q-2}2\int_{\Omega}\varphi^{2}(x)(v^n(x))^q|\n
  \eta(x)|^2\,dx +\int_{\Omega}\varphi^{2^*}(x)V(x)\eta^2(x)
  v^2(x)(v^n(x))^{q-2}\,dx.
\end{align}
Estimate (\ref{eq:64}) applied to $\eta(v^n)^{\frac q2-1}v$ gives 
\begin{align}
\label{eq:5bk}
\int_{\Omega} \varphi^{2^*}(x)|V(x)|\big[\eta(x)(v^n(x))^{\frac q2-1}v(x)\big]^2\,dx
&\leq
\frac{C(q)}{2}\int_{\Omega}\varphi^{2}(x)\big|\n(\eta(v^n)^{\frac q2-1}v)(x)\big|^2\,dx\\
&\quad\notag
+\ell_q\int_{\Omega}\varphi^{2^*}(x)(v^n(x))^{q-2}v^2(x)\eta^2(x)\,dx. 
\end{align}
Using  (\ref{eq:5bk}) to estimate the term with $V$ in (\ref{eq:22bk}), 
 (\ref{eq:63}) yields
\begin{align*}
  \bigg(\int_{\Omega}& \varphi^{2^*}(x)
|v^n(x)|^{(\frac q2-1)2^*}|v(x)|^{2^*}\eta^{2^*}(x)\,dx\bigg)^{\frac 2{2^*}}
  \leq
  \frac{2\ell_q{C}_{N,a}}{C(q)}\int_{\Omega} \varphi^{2^*}(x)
\eta^2(x)|v^n(x)|^{q-2}v^2(x)\,dx\\
&\quad +\frac{4{C}_{N,a}
 (2+C(q))}{C(q)}\int_{\Omega} \varphi^{2}(x)|v^n(x)|^{q-2}v^2(x)|\n \eta(x)|^2\,dx\\
  &\quad +{C}_{N,a}(q-2)\int_{\Omega} \varphi^{2}(x)|v^n(x)|^q |\n
  \eta(x)|^2\,dx.
\end{align*}
Letting $n\to\infty$ in the above inequality, (\ref{eq:66}) follows. 
\end{pf}

\noindent
\begin{pfn}{Theorem \ref{t:bk}} Let $u$ be a weak
  $H^1(\Omega)$-solution to \eqref{eq:67}.  It is easy to verify that
  $\varphi(x):=|x|^{\sigma}\psi_1(x/|x|)\in H^1(\Omega)$ satisfies
  (in a weak $H^1(\Omega)$-sense and in a classical sense in
  $\Omega\setminus\{0\}$)
$$
-\Delta \varphi(x)-\frac{a(x/|x|)}{|x|^2}\,\varphi(x)=0.
$$
Then $v:=\frac{u}{\varphi}\in H^1_{\varphi}(\Omega)$ turns out to be a
weak solution to \eqref{eq:65}. Let $R>0$ be such that 
$$
\Omega'\Subset
\Omega'+B(0,2R)\Subset \Omega.
$$
Using Lemma \ref{l:bk1} in $\Omega_1:=\Omega'+B(0,R(2-r_1))\Subset
\Omega'+B(0,2R)$, $r_1=1$, with $q=q_1=2^*$, we infer that $v\in
L^{\frac{(2^*)^2}2}(\varphi^{2^*}\!\!,\Omega_1)$ and the following
estimate holds
\begin{multline*}
  \|v\|_{L^{\frac{(2^*)^2}2}(\varphi^{2^*}\!\!,\Omega_1)}\leq
  \widetilde C^{\frac1{q_1}} \big(\!\mathop{\rm
    diam}\Omega\big)^{\frac{\sigma(2-2^*)}{q_1}}\bigg(\frac{8}{C(q_1)}\frac{4}{
    (R r_1)^2}+ \frac{4(q_1+2)}{(R r_1)^2}
  +\frac{2\ell_{q_1}}{C(q_1)}\bigg)^{\!\!\frac1{q_1}}
  \|v\|_{L^{2^*}(\varphi^{2^*}\!\!,\Omega)}.
\end{multline*}
Using again Lemma \ref{l:bk1} in $\Omega_2:=\Omega'+B(0,R(2-r_1-r_2))\Subset
\Omega_1$, $r_2=\frac14$, with $q=q_2=(2^*)^2/2$, we infer that $v\in
L^{\frac{(2^*)^3}4}(\varphi^{2^*}\!\!,\Omega_2)$ and 
\begin{gather*}
  \|v\|_{L^{\frac{(2^*)^3}4}(\varphi^{2^*}\!\!,\Omega_2)}\leq
  \widetilde C^{\frac1{q_2}} \big(\!\mathop{\rm
    diam}\Omega\big)^{\frac{\sigma(2-2^*)}{q_2}}\bigg(\frac{8}{C(q_2)}\frac{4}{
    (R r_2)^2}+ \frac{4(q_2+2)}{(R r_2)^2}
  +\frac{2\ell_{q_2}}{C(q_2)}\bigg)^{\!\!\frac1{q_2}}
  \|v\|_{L^{q_2}(\varphi^{2^*}\!\!,\Omega_1)}\\
  \leq \Big[\widetilde C\big(\!\mathop{\rm
    diam}\Omega\big)^{{\sigma(2-2^*)}}\Big]^{\frac1{q_1}+\frac1{q_2}}
  \bigg(\frac{8}{C(q_1)}\frac{4}{ (R r_1)^2}+ \frac{4(q_1+2)}{(R
    r_1)^2}
  +\frac{2\ell_{q_1}}{C(q_1)}\bigg)^{\!\!\frac1{q_1}}\times\\
  \times\bigg(\frac{8}{C(q_2)}\frac{4}{ (R r_2)^2}+ \frac{4(q_2+2)}{(R
    r_2)^2} +\frac{2\ell_{q_2}}{C(q_2)}\bigg)^{\!\!\frac1{q_2}}
  \|v\|_{L^{2^*}(\varphi^{2^*}\!\!,\Omega)}.
\end{gather*}
Setting, for any  $n\in\N$, $n\geq 1$,
$$
q_n=\frac12\bigg(\frac{2^*}2\bigg)^{\!n},\quad
\Omega_n:=\Omega'+B\bigg(0,R\bigg(2-\sum_{k=1}^n r_k\bigg)\bigg),\quad\text{and}
\quad
r_n=\frac1{n^2},
$$
and using iteratively Lemma \ref{l:bk1}, we obtain that, for any $n\in\N$, $n\geq 1$, 
\begin{multline}\label{eq:68}
  \|v\|_{L^{q_{n+1}}(\varphi^{2^*}\!\!,\Omega')}\leq 
 \|v\|_{L^{q_{n+1}}(\varphi^{2^*}\!\!,\Omega_n)}\\
\leq \|v\|_{L^{2^*}(\varphi^{2^*}\!\!,\Omega)}
   \Big[\widetilde C\big(\!\mathop{\rm
    diam}\Omega\big)^{{\sigma(2-2^*)}}\Big]^{\sum\limits_{k=1}^n\frac1{q_k}}
\prod_{k=1}^n  \bigg(\frac{8}{C(q_k)}\frac{4}{ (R r_k)^2}+ \frac{4(q_k+2)}{(R
    r_k)^2}
  +\frac{2\ell_{q_k}}{C(q_k)}\bigg)^{\!\!\frac1{q_k}}.
\end{multline}
We notice that
\begin{align}\label{eq:71}
\prod_{k=1}^n \bigg(\frac{8}{C(q_k)}\frac{4}{ (R r_k)^2}+ \frac{4(q_k+2)}{(R
    r_k)^2}
  +\frac{2\ell_{q_k}}{C(q_k)}\bigg)^{\!\!\frac1{q_k}}=\exp\bigg[\sum_{k=1}^n b_k\bigg]
\end{align}
where
$$
b_k=
\frac1{q_k}\log
\bigg(\frac{16 k^4}{R^2 C(q_k)}+ \frac{4k^4(q_k+2)}{R^2}
  +\frac{2\ell_{q_k}}{C(q_k)}\bigg),
$$
and, for some constant $C=C\big(N,a,\|V\|_{L^s(\varphi^{2^*}\!\!,\Omega)},\mathop{\rm
    dist}(\Omega',\partial\Omega)\big)>0$,
$$
b_k\sim
2\bigg(\frac{2}{2^*}\bigg)^{\!\!k}\log\left[C 
\bigg(\frac12\bigg(\frac{2^*}2\bigg)^{\!\!k}\bigg)^{\!\frac{2s}{2s-N}}\right]
\quad\text{as }k\to+\infty.
$$
Hence $\sum_{n=1}^\infty b_n$ converges to some positive sum depending
only on  $\|V\|_{L^s(\varphi^{2^*}\!\!,\Omega)}$, 
$\mathop{\rm dist}(\Omega',\partial\Omega)$, $N$, and $a$, hence
$$
\lim_{n\to+\infty}\Big[\widetilde C\big(\!\mathop{\rm
    diam}\Omega\big)^{{\sigma(2-2^*)}}\Big]^{\sum\limits_{k=1}^n\frac1{q_k}}
\prod_{k=1}^n  \bigg(\frac{8}{C(q_k)}\frac{4}{ (R r_k)^2}+ \frac{4(q_k+2)}{(R
    r_k)^2}
  +\frac{2\ell_{q_k}}{C(q_k)}\bigg)^{\!\!\frac1{q_k}}
$$
is finite and depends only on $N$, $a$,
$\|V\|_{L^s(\varphi^{2^*}\!\!,\Omega)}$, and $\mathop{\rm
  dist}(\Omega',\partial\Omega)$.
Hence, from \eqref{eq:68}, we deduce that there exists  a positive constant
  $C$
  depending only on $N$, $a$,
  $\|V\|_{L^s(\varphi^{2^*}\!\!,\Omega)}$, $\mathop{\rm dist}
  (\Omega',\partial\Omega)$, and $\mathop{\rm diam}\Omega$, such that
\begin{equation*}
  \|v\|_{L^{q_{n+1}}(\varphi^{2^*}\!\!,\Omega')}\leq  
C\,\|v\|_{L^{2^*}(\varphi^{2^*}\!\!,\Omega)}\quad\text{for all }n\in\N.
\end{equation*}
Letting $n\to+\infty$ we deduce that $v$ is essentially bounded in
$\Omega'$ with respect to the measure $\varphi^{2^*}\!dx$ and
\begin{equation*}
 \|v\|_{L^{\infty}(\varphi^{2^*}\!\!,\Omega')}\leq  
  C\,\|v\|_{L^{2^*}(\varphi^{2^*}\!\!,\Omega)}=C\,
\|u\|_{L^{2^*}(\Omega)},
\end{equation*}
where $\|v\|_{L^{\infty}(\varphi^{2^*}\!\!,\Omega')}$ denotes the
essential supremum of $v$ with respect to the measure
$\varphi^{2^*}\!dx$. Since $\varphi^{2^*}\!dx$ is absolutely
continuous with respect to the Lebesgue measure and viceversa, there holds
$\|v\|_{L^{\infty}(\varphi^{2^*}\!\!,\Omega')}=\|v\|_{L^{\infty}(\Omega')}$, hence 
 $v\in L^{\infty}(\Omega')$ and
\begin{equation*}
  \|v\|_{L^{\infty}(\Omega')}\leq  
  C\,
\|u\|_{L^{2^*}(\Omega)},
\end{equation*}
thus completing the proof.
\end{pfn}

\noindent
If the potential $V$ in equation \eqref{eq:67} belongs to
$L^{N/2}(\varphi^{2^*}\!\!,\Omega)$ (but to
$L^s(\varphi^{2^*}\!\!,\Omega)$ for no $s>N/2$), although we can no
more derive an $L^{\infty}$- bound for $u/\varphi$, we can obtain for  
$u/\varphi$ as high summability as we like.

\begin{Theorem}\label{t:bkn2}
  Let $\Omega$ be a bounded domain containing $0$, $a\in
  L^{\infty}({\mathbb S}^{N-1})$ satisfying $\Lambda_N(a)<1$, and
  $V\in L^{N/2}(\varphi^{2^*}\!\!,\Omega)$. Then, for any $\Omega'
  \Subset \Omega$ and for any weak $H^1(\Omega)$-solution $u$ to
  \eqref{eq:67}, there holds $\frac {u}{\varphi}\in
  L^{q}(\varphi^{2^*}\!\!,\Omega')$ for all $1\leq q<+\infty$.
\end{Theorem}
\begin{pf}
  The proof follows closely  the proofs of Theorem \ref{t:bk} and Lemma
  \ref{l:bk1}. However, since we only require $V\in
  L^{N/2}(\varphi^{2^*}\!\!,\Omega)$, we have that for any $q$ there
  exists $\ell_q$ such that 
\begin{align*}
  \int\limits_{|V(x)|\geq \ell_q} \varphi^{2^*}(x)|V(x)|^{\frac
    N2}\,dx \leq\min\bigg\{\frac{C_{N,a}^{-1}}{8},
  \frac{2C_{N,a}^{-1}}{q+4}\bigg\}^{\!\!\frac N2}, 
\end{align*}
but we can no more estimate $\ell_q$ in terms of $q$, as we did in
\eqref{eq:70} thanks to the summability assumption
$V\in L^s(\varphi^{2^*}\!\!,\Omega)$ for some $s>N/2$. Hence we still
arrive at an estimate of type \eqref{eq:68} but we have no control on
the product in \eqref{eq:71} as $n\to+\infty$.
\end{pf}

\section{Behavior of solutions at singularities}\label{sec:bound}

The procedure followed in this section to prove Theorem \ref{t:1}
relies in comparison methods and separation of variables. Indeed we
will evaluate the asymptotics of solutions to problem (\ref{eq:45}) by
trapping them between functions which solve analogous problems with
radial perturbing potentials. To this aim, the first step consists in
deriving the asymptotic behavior of solutions to Schr\"odinger
equations with a potential which is given by a radial perturbation of
the dipole-type singular term. In this case, it is possible to expand
the solution in Fourier series, thus separating the radial and angular
variables, and to estimate the behavior of the Fourier coefficients in
order to establish which of them is dominant near the singularity.

\begin{Proposition}\label{p:1}
  Let $a\in L^{\infty}({\mathbb S}^{N-1})$ be such that
  $\Lambda_N(a)<1$, $R>0$, and $u\in H^1(B(0,R))$, $u\geq 0$ a.e. in
  $B(0,R)$, $u\not\equiv 0$, be a weak $H^1$-solution to
\begin{equation}\label{eq:1}
  -\Delta u(x)=\bigg[\frac{a(x/|x|)}
  {|x|^2}+h(|x|)\bigg]\,u(x)\quad\text{ in }B(0,R),
\end{equation} 
where $h\in  L^{\infty}_{\rm loc}(0,R)\cap L^p(0,R)$ for
some $p>N/2$. 
Then, for any $r\in(0,R)$, there exists a positive constant $C$ (depending
on $h$, $R$, $r$, $a$, $\e$, and $u$) such that
\begin{align*}
\frac1C|x|^{\sigma}\leq u(x)\leq
C|x|^{\sigma}\quad\text{ for all }x\in B(0,r)\setminus\{0\},
\end{align*}
where $\sigma$ is defined in (\ref{eq:sigma}).  Moreover, there exists a 
positive constant $\widetilde C$ (depending on $u$, $h$, $N$, and $a$)
such that,
for any
$\theta\in{{\mathbb S}^{N-1}}$, 
\begin{equation}\label{eq:40}
\lim_{\rho\to0^+}u(\rho\,\theta)
\rho^{-\sigma}=\widetilde C
\,\psi_1(\theta),
\end{equation}
and, for any $r\in(0,R)$,
\begin{multline}\label{eq:58}
\widetilde C=\int_{{\mathbb
      S}^{N-1}} \bigg(r^{-\sigma} u(r\eta)+\int_0^r 
{{\frac{s^{1-\sigma}}
  {2\sigma+N-2}}}\,h(s)u(s\,\eta)\,ds
\\
-r^{-2\sigma-N+2}\int_0^r 
{{\frac{s^{N-1+\sigma}}
  {2\sigma+N-2}}}\,h(s)u(s\,\eta)\,ds
\bigg)\psi_1(\eta)
  \,dV(\eta).
\end{multline}
Furthermore, for any $r\in(0,R)$,  there exists a positive constant
$\bar C$ (depending on $h$, $R$, $r$, $a$, $\e$, but not on $u$) such
that
\begin{align}\label{eq:50}
  u(\rho\,\theta)\leq \bar
  C\,\|u\|_{H^1(B(0,R))}\rho^{\sigma},\quad\text{for all }0<\rho<r.
\end{align}
\end{Proposition}

\begin{pf}
  Let $r\in(0,R)$. We can assume, without loss of generality, that
  $R>1$ and $r=1$. Indeed,  setting $w(x):=u(rx)$, we
  notice that $w\in H^1(B(0,R/r))$ and weakly solves
$$
  -\Delta w(x)=\bigg[\frac{a(x/|x|)}
  {|x|^2}+\tilde h(|x|)\bigg]\,w(x)\quad\text{ in }B(0,R/r),
$$
where $\tilde h(\rho):=r^2 h(r\rho)$ satisfies $\tilde h\in
L^{\infty}_{\rm loc}(0,R/r)\cap L^p(0,R/r)$. Hence,
it is enough to prove the statement for $R>1$ and $r=1$, being the general 
case easily obtainable from scaling.

Let $R>1$, $r=1$ and $u\in H^1(B(0,R))$, $u\geq 0$ a.e. in $B(0,R)$,
$u\not\equiv 0$, be a weak solution of (\ref{eq:1}).  By standard
regularity theory, $u\in C^0\big(\overline{B(0,1)}\setminus
B(0,s)\big)$ for any $s\in(0,1)$.  For any $k\in\N\setminus\{0\}$,
let $\psi_k$   be a $L^2$-normalized
eigenfunction of the operator $-\D_{\mathbb
  S^{N-1}}-a(\theta)$ on the sphere associated to
the $k$-th eigenvalue $\mu_k$, i.e. satisfying \eqref{eq:2rad}.
We can choose the functions $\psi_k$ in such a way that they form an
orthonormal basis of $L^2({\mathbb S}^{N-1})$, hence $u$ can be expanded as 
\begin{equation}\label{sviluppo}
u(x)=u(\rho\,\theta)=\sum_{k=1}^\infty\varphi_k(\rho)\psi_k(\theta),
\end{equation}
where $\rho=|x|\in(0,1]$, $\theta=x/|x|\in{{\mathbb S}^{N-1}}$, and
\begin{equation}\label{eq:3}
\varphi_k(\rho)=\int_{{\mathbb S}^{N-1}}u(\rho\,\theta)\psi_k(\theta)\,dV(\theta).
\end{equation}
The Parseval identity yields
$$
\int_{{\mathbb
    S}^{N-1}}|u(\rho\,\theta)|^2\,dV(\theta)=
\sum_{k=1}^{\infty}|\varphi_k(\rho)|^2,\quad\text{for all }0<\rho\leq 1,
$$
and hence
\begin{equation}\label{eq:parseval}
  \|u\|^2_{L^2(B(0,1))}=  \int_0^1\rho^{N-1}\bigg(\sum_{k=1}^{\infty}|\varphi_k(\rho)|^2
  \bigg)\,d\rho= \sum_{k=1}^{\infty}\int_0^1\rho^{N-1}|\varphi_k(\rho)|^2
  \,d\rho.
\end{equation}
Equations \eqref{eq:1} and \eqref{eq:2rad} imply that, for every $k$,
\begin{equation*}
\varphi_k''(\rho)+\frac{N-1}{\rho}\varphi_k^\prime(\rho)-
\frac{\mu_k}{\rho^2}\varphi_k(\rho)=h(\rho)\varphi_k(\rho),\quad\text{in }(0,1).
\end{equation*}
A direct calculation shows that, for some $c_1^k,c_2^k\in\R$,
\begin{equation}\label{eq:42}
\varphi_k(\rho)=\rho^{\sigma^+_k}
\bigg(c_1^k+\int_\rho^1\frac{s^{-\sigma^+_k+1}}{\sigma^+_k-\sigma^-_k}
h(s)\varphi_k(s)\,ds\bigg)+\rho^{\sigma^-_k}
\bigg(c_2^k+\int_\rho^1\frac{s^{-\sigma^-_k+1}}{\sigma^-_k-\sigma^+_k}
h(s)\varphi_k(s)\,ds\bigg),
\end{equation}
where 
\begin{equation}\label{eq:56}
\sigma^+_k=-\frac{N-2}{2}+\sqrt{\bigg(\frac{N-2}
  {2}\bigg)^{\!\!2}+\mu_k}\quad\text{and}\quad
\sigma^-_k=-\frac{N-2}{2}-\sqrt{\bigg(\frac{N-2}{2}\bigg)^{\!\!2}+\mu_k}.
\end{equation}
For the sake of notation, we set
\begin{equation}\label{eq:14}
A_k(\rho)=\rho^{\sigma^+_k}
\int_\rho^1\frac{s^{-\sigma^+_k+1}}{\sigma^+_k-\sigma^-_k}
\,h(s)\varphi_k(s)\,ds
\end{equation}
and 
\begin{equation*}
B_k(\rho)=\rho^{\sigma^-_k}
\bigg(c_2^k+\int_\rho^1\frac{s^{-\sigma^-_k+1}}{\sigma^-_k-\sigma^+_k}
\,h(s)\varphi_k(s)\,ds\bigg),
\end{equation*}
so that
\begin{equation}\label{eq:11}
\varphi_k(\rho)=c_1^k\rho^{\sigma^+_k}+A_k(\rho)+B_k(\rho).
\end{equation}
Without loss of generality, we can assume that
$$\frac{N}{2}<p<\frac{N}{2-\frac{1}{3}\sqrt{\big(\frac{N-2}
{2}\big)^2+\mu_1}},$$
so, setting $\e=2-\frac{N}{p}$, $0<\e<\frac{1}{3}(\sigma_k^++\frac{N-2}
{2})$, for every $k$. 
From H\"older's
inequality and Lemma~\ref{l:stiautf}, it follows that
\begin{gather}\label{eq:6}
  |A_k(\rho)|=\rho^{\sigma^+_k}\bigg|
  \int_\rho^1\frac{s^{-\sigma^+_k+1}}{\sigma^+_k-\sigma^-_k}
  h(s)\bigg(\int_{{\mathbb
      S}^{N-1}}u(s\,\theta)\psi_k(\theta)\,dV(\theta)\bigg)\,ds\bigg|\\
  \notag\leq \frac{ 
    C_1\rho^{\sigma^+_k}|\mu_k|^{\big\lfloor\!\frac{N-1}4\!\big\rfloor+1}}{\sigma^+_k-\sigma^-_k}
  \int_\rho^1s^{-\sigma^+_k+1-\frac{N-1}{p}-\frac{N-1}{2^*}}
  \bigg(\int_{{\mathbb
      S}^{N-1}}s^{\frac{N-1}{p}+\frac{N-1}{2^*}}|h(s)|
  |u(s\,\theta)|\,dV(\theta)\bigg)\,ds\\
  \notag\leq\frac{\omega_N^{\frac{N+2}{2N}}C_1\rho^{\sigma^+_k}|\mu_k|^{\big\lfloor\!
      \frac{N-1}4\!\big\rfloor+1}}{\sigma^+_k-\sigma^-_k}\,
  \frac{\|h\|_{L^{p}(B(0,1))}}{\omega_N}\,\|u\|_{L^{2^*}(B(0,1))}
  \bigg(\int_\rho^1\! s^{(-\sigma^+_k+1-\frac{N-1}{p}-\frac{N-1}{2^*})
    \frac{2^*p}{2^*p-2^*-p}}\,ds\bigg)^{\!\!1-\frac1p-\frac1{2^*}}\\
  \notag=\frac{\omega_N^{\frac{N+2}{2N}}C_1\rho^{\sigma^+_k}|\mu_k|^{\big\lfloor\!
      \frac{N-1}4\!\big\rfloor+1}}{\sigma^+_k-\sigma^-_k}\,
  \frac{\|h\|_{L^{p}(B(0,1))}}{\omega_N}\,\|u\|_{L^{2^*}(B(0,1))}
  \Bigg[\frac{\rho^{\frac{2^*p}{2^*p-2^*-p}(-\sigma_k^+-\frac{N-2}2+2-\frac{N}p)}}
  {\frac{2^*p}{2^*p-2^*-p}\big(\frac{N}p-2+\sigma_k^+
    +\frac{N-2}2\big)}\Bigg]^{1-\frac1p-\frac1{2^*}}\\
  \notag\leq\frac{\omega_N^{\frac{N+2}{2N}}C_1
    |\mu_k|^{\big\lfloor\!\frac{N-1}4\!\big\rfloor+1}}{\sigma^+_k-\sigma^-_k}
  \,\frac{\|h\|_{L^{p}(B(0,1))}}{\omega_N}\, \|u\|_{L^{2^*}(B(0,1))}
  \frac{\rho^{\e-\frac{N-2}2}}
  {\Big[\frac{2^*p}{2^*p-2^*-p}\big(\sigma_k^+
    +\frac{N-2}2-\e\big)\Big]^{1-\frac1p-\frac1{2^*}}},
\end{gather}
where $\omega_N=\int_{{\mathbb
S}^{N-1}}dV(\theta)$.
In particular
\begin{equation}\label{eq:5}
A_k(\rho)=o(\rho^{\sigma_k^-})\quad\text{as }\rho\to0^+.
\end{equation}
Moreover
\begin{align}\label{eq:9}
  &\int_0^1\big|{s^{-\sigma^-_k+1}}h(s)\varphi_k(s)\big|\,ds\\
  &\notag\!\!\!\!\!\!\leq
  C_1|\mu_k|^{\big\lfloor\!\frac{N-1}4\!\big\rfloor+1}
  \int_0^1s^{-\sigma^-_k+1-\frac{N-1}{p}-\frac{N-1}{2^*}}
  \bigg(\int_{{\mathbb
      S}^{N-1}}s^{\frac{N-1}{p}+\frac{N-1}{2^*}}|h(s)||u(s\,\theta)|\,dV(\theta)\bigg)\,ds\\
  &\notag\!\!\!\!\!\!\leq
  \omega_N^{\frac{N+2}{2N}}C_1|\mu_k|^{\big\lfloor\!\frac{N-1}4\!\big\rfloor+1}
  \,\frac{\|h\|_{L^{p}(B(0,1))}}{\omega_N}\,
  \|u\|_{L^{2^*}(B(0,1))}\bigg[\int_0^1
  \!s^{(-\sigma^-_k+1-\frac{N-1}{p}-\frac{N-1}{2^*})
    \frac{2^*p}{2^*p-2^*-p}}\,ds\bigg]^{1-\frac1p-\frac1{2^*}}\\
  &\notag\!\!\!\!\!\!=\omega_N^{\frac{N+2}{2N}}C_1
  |\mu_k|^{\big\lfloor\!\frac{N-1}4\!\big\rfloor+1}\,\frac{\|h\|_{L^{p}(B(0,1))}}{\omega_N}
  \,\|u\|_{L^{2^*}(B(0,1))} \frac{1}{\big[\frac{2^*p}{2^*p-2^*-p}
   \big(\e-\sigma_k^--\frac{N-2}2\big)\big]^{1-\frac1p-\frac1{2^*}}} <\infty
\end{align}
due to inequality $\sigma^-_k<-\frac{N-2}{2}$. As a consequence, 
\begin{equation}\label{eq:7}
\lim_{\rho\to 0^+}\int_\rho^1\frac{s^{-\sigma^-_k+1}}{\sigma^-_k-\sigma^+_k}
\, h(s)\varphi_k(s)\,ds
\end{equation}
is finite. Since $u\in L^{2^*}(B(0,1))$, from (\ref{eq:5}),
(\ref{eq:7}), and the fact that $\rho^{\sigma^-_k}\psi_k(\theta)\not
\in L^{2^*}(B(0,1))$, we conclude that there must be
\begin{equation}\label{eq:43}
c_2^k=-\int_0^1\frac{s^{-\sigma^-_k+1}}{\sigma^-_k-\sigma^+_k}
\,h(s)\varphi_k(s)\,ds,
\end{equation}
hence 
\begin{align}\label{eq:8}
B_k(\rho)=\rho^{\sigma^-_k}
\int_0^\rho\frac{s^{-\sigma^-_k+1}}{\sigma^+_k-\sigma^-_k}
h(s)\varphi_k(s)\,ds.
\end{align}
Since $u\in
C^0\big(\overline{B(0,1)}\setminus B(0,s)\big)$ for any $s\in(0,1)$,
 it makes sense to evaluate $\varphi_k$ at $\rho=1$ and, 
from (\ref{eq:42}) and (\ref{eq:43}), we have that
\begin{equation}\label{eq:59}
c_1^k=\varphi_k(1)+\int_0^1\frac{s^{-\sigma^-_k+1}}{\sigma^-_k-\sigma^+_k}
\,h(s)\varphi_k(s)\,ds.
\end{equation}
From (\ref{eq:42}), (\ref{eq:8}), and (\ref{eq:59}), we deduce that 
\begin{align}\label{eq:60}
\varphi_k(\rho)=\rho^{\sigma^+_k}
\bigg(\varphi_k(1)&+\int_0^1\frac{s^{-\sigma^-_k+1}}{\sigma^-_k-\sigma^+_k}
\,h(s)\varphi_k(s)\,ds
+\int_\rho^1\frac{s^{-\sigma^+_k+1}}{\sigma^+_k-\sigma^-_k}
h(s)\varphi_k(s)\,ds\bigg)\\
\notag
&+\rho^{\sigma^-_k}
\int_0^\rho\frac{s^{-\sigma^-_k+1}}{\sigma^+_k-\sigma^-_k}
h(s)\varphi_k(s)\,ds.
\end{align}
From above, (\ref{eq:9}), and standard elliptic estimates (which allow
to estimate $u$ outside the singularity in terms of its $H^1$-norm) we
obtain that, for some positive constant $\tilde c$ depending only on
$N$, $R$, $a$, and~$h$, 
\begin{equation}\label{eq:46}
|c_1^k|\leq \tilde c\,
|\mu_k|^{\big\lfloor\!\frac{N-1}4\!\big\rfloor+1}\|u\|_{H^1(B(0,R))}
\Big[1+\Big[{\textstyle\frac{2^*p}{2^*p-2^*-p}(\e-\sigma_k^--\frac{N-2}2)
\Big]^{-1+\frac1p+\frac1{2^*}}}\Big].
\end{equation}
Arguing as in (\ref{eq:9}), we find that 
\begin{equation}\label{eq:10}
|B_k(\rho)|\leq\frac{\omega_N^{\frac{N+2}{2N}}
C_1|\mu_k|^{\big\lfloor\!\frac{N-1}4\!\big\rfloor+1}}{\sigma^+_k-\sigma^-_k}
\,\frac{\|h\|_{L^{p}(B(0,1))}}{\omega_N}\,\|u\|_{L^{2^*}(B(0,1))} 
\frac{\rho^{\e-\frac{N-2}2}}
{\big[\frac{2^*p}{2^*p-2^*-p}(\e-\sigma_k^--\frac{N-2}2)\big]^{1-\frac1p-\frac1{2^*}}}.
\end{equation}
From (\ref{eq:11}), (\ref{eq:6}),
and (\ref{eq:10}), we can estimate $\varphi_k$ as
\begin{equation}\label{eq:12}
|\varphi_k(\rho)|\leq|c_1^k|\rho^{\sigma_k^+}+
\frac{\a_k|\mu_k|^{\big\lfloor\!\frac{N-1}4\!\big\rfloor
+1}}{\sigma_k^+-\sigma_k^-}\,\frac{\|h\|_{L^{p}(B(0,1))}}{\omega_N}\,\|u\|_{L^{2^*}(B(0,1))}
\rho^{\e-\frac{N-2}2}
\end{equation}
where
$$
\a_k:=\frac{\omega_N^{\frac{N+2}{2N}}C_1}{\big[\frac{2^*p}{2^*p-2^*-p}(\sigma_k^+
+\frac{N-2}2-\e)\big]^{1-\frac1p-\frac1{2^*}}}+
\frac{\omega_N^{\frac{N+2}{2N}}C_1}{\big[\frac{2^*p}{2^*p-2^*-p}
(\e-\sigma_k^--\frac{N-2}2)\big]^{1-\frac1p-\frac1{2^*}}}.
$$
\medskip\noindent
{\bf Claim 1:} there holds
\begin{align}\label{eq:13}
&|A_k(\rho)+B_k(\rho)|\leq
|c_1^k|\rho^{\sigma_k^+}\sum_{i=1}^{j_k-1}\bigg(\frac{\|h\|_{L^{p}(B(0,1))}}
{\omega_N(\sigma_k^+-\sigma_k^-)}\frac{2}{\e^{\frac{p-1}{p}}}\bigg)^{\!i}\\
&\notag\quad+\a_k|\mu_k|^{\big\lfloor\!\frac{N-1}4\!\big\rfloor
+1}\|u\|_{L^{2^*}(B(0,1))}\bigg(\frac{\|h\|_{L^{p}(B(0,1))}}
{\omega_N(\sigma_k^+-\sigma_k^-)}\bigg)^{j_k}
\rho^{j_k\e-\frac{N-2}2}\prod_{i=2}^{j_k}
\frac2{\Big[\sqrt{\big(\frac{N-2}{2}\big)^{\!2}+\mu_k}-i\e\Big]^{\frac{p-1}{p}}},
\end{align}
where 
$$
j_k:= \bigg\lfloor \frac1\e\sqrt{\Big(\frac{N-2}{2}
\Big)^2+\mu_k}\bigg\rfloor-1,
$$
i.e. $j_k$ is the unique integer number such that
$$
\sqrt{\Big(\frac{N-2}{2}
\Big)^2+\mu_k}-2\e<j_k\e\leq \sqrt{\Big(\frac{N-2}{2}
\Big)^2+\mu_k}-\e.
$$
Notice that $\e\leq\frac13\big(\sigma_k^++\frac{N-2}2\big)$ implies
$j_k\geq2$. 

To prove the claim, we observe that, from (\ref{eq:14})
and (\ref{eq:12}) it follows that
\begin{align}\label{eq:16}
&|A_k(\rho)|\leq|c_1^k|\rho^{\sigma_k^+}\frac{\|h\|_{L^{p}(B(0,1))}}{\omega_N(\sigma_k^+-\sigma_k^-)}
\Big(\int_\rho^1s^{(1-\frac{N-1}{p})\frac{p}{p-1}}\,ds\Big)^{\frac{p-1}{p}}\\
\notag&+\a_k|\mu_k|^{\big\lfloor\!\frac{N-1}4\!\big\rfloor
+1}\|u\|_{L^{2^*}(B(0,1))}
\bigg(\frac{\|h\|_{L^{p}(B(0,1))}}{\omega_N(\sigma_k^+-\sigma_k^-)}\bigg)^2
\rho^{\sigma_k^+}\Big(\int_\rho^1
s^{(-\sigma_k^++1+\e-\frac{N-2}2-\frac{N-1}{p})\frac{p}{p-1}}\,ds\Big)^{\frac{p-1}{p}}\\
&\notag\leq|c_1^k|\rho^{\sigma_k^+}
\frac{\|h\|_{L^{p}(B(0,1))}}{\omega_N(\sigma_k^+-\sigma_k^-)}\frac{1}{(\e\frac{p}{p-1})^{\frac{p-1}{p}}}\!\\
&\notag+
\a_k|\mu_k|^{\big\lfloor\!\frac{N-1}4\!\big\rfloor
+1}\|u\|_{L^{2^*}(B(0,1))}
\bigg(\frac{\|h\|_{L^{p}(B(0,1))}}{\omega_N(\sigma_k^+-\sigma_k^-)}\bigg)^2
\frac{\rho^{2\e-\frac{N-2}2}}{\Big[\big(\sigma_k^+-2\e+\frac{N-2}2\big)
\frac{p}{p-1}\Big]^{\frac{p-1}{p}}}.
\end{align}
In a similar way, from (\ref{eq:8}) and (\ref{eq:12}) we deduce that
\begin{align}\label{eq:17}
|B_k(\rho)|\leq &\,|c_1^k|\rho^{\sigma_k^+}
\frac{\|h\|_{L^{p}(B(0,1))}}{\omega_N(\sigma_k^+-\sigma_k^-)}
\frac{1}{\big[(\e+\sigma_k^+-\sigma_k^-)\frac{p}{p-1}\big]^{\frac{p-1}{p}}}\\
\notag&\quad+\a_k|\mu_k|^{\big\lfloor\!\frac{N-1}4\!\big\rfloor
+1}\|u\|_{L^{2^*}(B(0,1))}
\bigg(\frac{\|h\|_{L^{p}(B(0,1))}}{\omega_N(\sigma_k^+-\sigma_k^-)}\bigg)^2
\frac{\rho^{2\e-\frac{N-2}2}}
{\Big[\big(2\e-\sigma_k^--\frac{N-2}2\big)\frac{p}{p-1}\Big]^{\frac{p-1}{p}}}.
\end{align}
Summing up (\ref{eq:16}) and (\ref{eq:17}), we obtain
\begin{align*}
|A_k(\rho)+B_k(\rho)|&\leq
|c_1^k|\rho^{\sigma_k^+}
\frac{\|h\|_{L^{p}(B(0,1))}}{\omega_N(\sigma_k^+-\sigma_k^-)}
\frac{2}{\e^{\frac{p-1}{p}}}\\
&\quad+\a_k|\mu_k|^{\big\lfloor\!\frac{N-1}4\!\big\rfloor
+1}\|u\|_{L^{2^*}(B(0,1))}
\bigg(\frac{\|h\|_{L^{p}(B(0,1))}}{\omega_N(\sigma_k^+-\sigma_k^-)}\bigg)^2
\frac{2\rho^{2\e-\frac{N-2}2}}
{\big(\sigma_k^++\frac{N-2}2-2\e\big)^{\frac{p-1}{p}}},
\end{align*}
and hence, from (\ref{eq:11}),
\begin{align}\label{eq:18}
|\varphi_k(\rho)|&\leq |c_1^k|\rho^{\sigma_k^+}\bigg(1+\frac{\|h\|_{L^{p}(B(0,1))}}{\omega_N(\sigma_k^+-\sigma_k^-)}
\frac{2}{\e^{\frac{p-1}{p}}}\bigg)+\\
\notag&\quad\a_k|\mu_k|^{\big\lfloor\!\frac{N-1}4\!\big\rfloor
+1}\|u\|_{L^{2^*}(B(0,1))}
\bigg(\frac{\|h\|_{L^{p}(B(0,1))}}{\omega_N(\sigma_k^+-\sigma_k^-)}\bigg)^2
\frac{2\rho^{2\e-\frac{N-2}2}}
{\big(\sqrt{\big(\frac{N-2}2\big)^{\!2}+\mu_k}-2\e\big)^{\frac{p-1}{p}}}.
\end{align}
Using (\ref{eq:18}), we can improve our estimates of $A_k(\rho)$ and
$B_k(\rho)$ thus obtaining
\begin{align*}
|A_k(\rho)|&\leq |c_1^k|\rho^{\sigma_k^+}
\left[\frac{\|h\|_{L^{p}(B(0,1))}}{\omega_N(\sigma_k^+-\sigma_k^-)}
\frac{1}{\e^{\frac{p-1}{p}}}+
\bigg(\frac{\|h\|_{L^{p}(B(0,1))}}{\omega_N(\sigma_k^+-\sigma_k^-)}\bigg)^2
\frac{2}{\e^{2\frac{p-1}{p}}}\right]\\
&\quad+\a_k|\mu_k|^{\big\lfloor\!\frac{N-1}4\!\big\rfloor
+1}\|u\|_{L^{2^*}(B(0,1))}
\bigg(\frac{\|h\|_{L^{p}(B(0,1))}}{\omega_N(\sigma_k^+-\sigma_k^-)}\bigg)^3
\times\\
&\qquad\qquad\times \frac{\rho^{3\e-\frac{N-2}2}}{\Big[\sqrt{\big(
\frac{N-2}2\big)^{\!2}+\mu_k}-3\e\Big]^{\frac{p-1}{p}}}\frac{2}
{\Big[\sqrt{\big(\frac{N-2}2\big)^{\!2}+\mu_k}-2\e\Big]^{\frac{p-1}{p}}}
\end{align*}
and
\begin{align*}
|B_k(\rho)|&\leq |c_1^k|\rho^{\sigma_k^+}
\left[\frac{\|h\|_{L^{p}(B(0,1))}}{\omega_N(\sigma_k^+-\sigma_k^-)}
\frac{1}{\e^{\frac{p-1}{p}}}+
\bigg(\frac{\|h\|_{L^{p}(B(0,1))}}{\omega_N(\sigma_k^+-\sigma_k^-)}\bigg)^2
\frac{2}{\e^{2\frac{p-1}{p}}}\right]\\
&\quad+\a_k|\mu_k|^{\big\lfloor\!\frac{N-1}4\!\big\rfloor
+1}\|u\|_{L^{2^*}(B(0,1))}
\bigg(\frac{\|h\|_{L^{p}(B(0,1))}}{\omega_N(\sigma_k^+-\sigma_k^-)}\bigg)^3
\times\\
&\qquad\qquad\times \frac{\rho^{3\e-\frac{N-2}2}}{\Big[\sqrt{\big(
\frac{N-2}2\big)^{\!2}+\mu_k}+3\e\Big]^{\frac{p-1}{p}}}\frac{2}
{\Big[\sqrt{\big(\frac{N-2}2\big)^{\!2}+\mu_k}-2\e\Big]^{\frac{p-1}{p}}}.
\end{align*}
Summing up we find that
\begin{align*}
|A_k(\rho)+B_k(\rho)|&\leq |c_1^k|\rho^{\sigma_k^+}
\left[\frac{\|h\|_{L^{p}(B(0,1))}}{\omega_N(\sigma_k^+-\sigma_k^-)}
\frac{2}{\e^{\frac{p-1}{p}}}+
\bigg(\frac{\|h\|_{L^{p}(B(0,1))}}{\omega_N(\sigma_k^+-\sigma_k^-)}\bigg)^2
\frac{4}{\e^{2\frac{p-1}{p}}}\right]\\
&\quad+\a_k|\mu_k|^{\big\lfloor\!\frac{N-1}4\!\big\rfloor
+1}\|u\|_{L^{2^*}(B(0,1))}
\bigg(\frac{\|h\|_{L^{p}(B(0,1))}}{\omega_N(\sigma_k^+-\sigma_k^-)}\bigg)^3\times\\
&\qquad\qquad\times
\frac{4\rho^{3\e-\frac{N-2}2}}{\Big[\sqrt{\big(
\frac{N-2}2\big)^{\!2}+\mu_k}-3\e\Big]^{\frac{p-1}{p}}
\Big[\sqrt{\big(\frac{N-2}2\big)^{\!2}+\mu_k}-2\e\Big]^{\frac{p-1}{p}}}.
\end{align*}
An iteration of the above argument $(j_k-1)$ times easily leads to estimate
(\ref{eq:13}). Claim 1 is thereby proved.

\medskip\noindent {\bf Claim 2:} the function $s \mapsto
s^{-\sigma^+_k+1}h(s)\varphi_k(s)$
belongs to $L^1(0,1)$ and 
\begin{equation}\label{eq:20}
\lim_{\rho\to0^+}\rho^{-\sigma_k^+}\varphi_k(\rho)=
c_1^k+ \int_0^1\frac{s^{-\sigma^+_k+1}}{\sigma^+_k-\sigma^-_k} h(s)\varphi_k(s)\,ds.
\end{equation}
Indeed, from (\ref{eq:13}),
(\ref{eq:11}), (\ref{eq:46}), and the choice of $j_k$, it follows that
\begin{equation}\label{eq:21}
|\varphi_k(\rho)|\leq d_k\|u\|_{H^1(B(0,1))} \rho^{j_k\e-\frac{N-2}2},
\end{equation}
for some positive constant $d_k$ depending on $k$ (and on $a$, $R$, $h$
$N$). We distinguish now two cases. If $j_k\e<
\sqrt{\big(\frac{N-2}{2} \big)^{\!2}+\mu_k}-\e$, then from (\ref{eq:14}),
(\ref{eq:8}), and (\ref{eq:21}) we derive that
\begin{equation*}
|A_k(\rho)+B_k(\rho)|\leq d_k'\|u\|_{H^1(B(0,1))}\rho^{(j_k+1)\e-\frac{N-2}2},
\end{equation*}
for some other positive constant $d_k'$ depending on $k$ (and on $a$, $R$, $h$,
$N$), and hence, by (\ref{eq:46}) and the choice of $j_k$,
\begin{equation}\label{eq:23}
|\varphi_k(\rho)|\leq d_k''\|u\|_{H^1(B(0,1))}\rho^{(j_k+1)\e-\frac{N-2}2}.
\end{equation}
Estimate (\ref{eq:23}) and the choice of $j_k$ imply that the function
$s \mapsto s^{-\sigma^+_k+1}h(s)\varphi_k(s)$ belongs to $L^1(0,1)$.
Moreover, from (\ref{eq:8}) and (\ref{eq:23}) it follows that
\begin{equation}\label{eq:47}
|B_k(\rho)|\leq d_k'''\|u\|_{H^1(B(0,1))}\rho^{(j_k+2)\e-\frac{N-2}2}=
o(\rho^{\sigma_k^+})\quad\text{as }\rho\to0^+,
\end{equation}
and the claim is proved.  If $j_k\e=\sqrt{\big(\frac{N-2}{2}
  \big)^{\!2}+\mu_k}-\e$, then from (\ref{eq:14}), (\ref{eq:8}), and
(\ref{eq:21}) we derive that
\begin{equation*}
|A_k(\rho)+B_k(\rho)|\leq \gamma_k\|u\|_{H^1(B(0,1))}
\rho^{\sigma_k^+}|\log\rho|^{\frac{p-1}{p}},
\end{equation*}
for some other positive constant $\gamma_k$ depending on $k$ (and on $a$, $R$, $h$,
$N$), and hence, 
\begin{equation}\label{eq:24}
|\varphi_k(\rho)|\leq \gamma_k'\|u\|_{H^1(B(0,1))}\rho^{\sigma_k^+}|\log\rho|^{\frac{p-1}{p}}.
\end{equation}
Estimate (\ref{eq:24})  implies that the function
$s \mapsto s^{-\sigma^+_k+1}h(s)\varphi_k(s)$ belongs to $L^1(0,1)$.
Moreover, from (\ref{eq:8}) and (\ref{eq:24}) it follows that
\begin{equation}\label{eq:48}
|B_k(\rho)|\leq\gamma_k''\|u\|_{H^1(B(0,1))}\rho^{\sigma_k^-}
\Big(\int_0^{\rho}s^{(\sigma_k^+-\sigma_k^-+\e)\frac{p}{p-1}-1}|\log s|\,ds\Big)^{\frac{p-1}{p}}=
o(\rho^{\sigma_k^+})\quad\text{as }\rho\to0^+,
\end{equation}
and claim 2 is proved also in this case. 

\medskip\noindent
Let us fix $\bar k$ such that
\begin{align}\label{eq:25}
\frac{\|h\|_{L^{p}(B(0,1))}}
{\omega_N(\sigma_k^+-\sigma_k^-)}\frac{2}{\e^{\frac{p-1}{p}}}
=\frac{\|h\|_{L^{p}(B(0,1))}}
{\omega_N\e^{\frac{p-1}{p}}\sqrt{\big(\frac{N-2}{2}
\big)^{\!2}+\mu_k}}<\frac13,
\quad\sigma_k^+>4\e\quad\text{and}\quad
\frac{\sigma_{k}^+}2>\sigma_1^+\,,\quad\forall\,k\geq\bar k.
\end{align}
From (\ref{eq:13}) and (\ref{eq:25}), it follows that, for all
$k\geq\bar k$ and  for some positive constant $C_2$ (depending only
on $N$, $h$, and $a$),
\begin{multline*}
|A_k(\rho)+B_k(\rho)|\leq
|c_1^k|\rho^{\sigma_k^+}\sum_{i=1}^{\infty}\Big(\frac13\Big)^i\\
+C_2\|u\|_{H^1(B(0,1))}\bigg(\frac{\|h\|_{L^{p}(B(0,1))}}
{\omega_N(\sigma_k^+-\sigma_k^-)}\bigg)^{\!j_k}
|\mu_k|^{\big\lfloor\!\frac{N-1}4\!\big\rfloor+1}
\frac{\rho^{j_k\e-\frac{N-2}2}\big(\frac2\e\big)^{j_k-1}}
{\prod_{i=2}^{j_k}\bigg(
\Big\lfloor{\e}^{-1}{\sqrt{\big(\frac{N-2}{2}\big)^{\!2}+
\mu_k}}\Big\rfloor-i\bigg)},
\end{multline*}
which yields, for all $k\geq\bar k$,
\begin{align}\label{eq:29}
|A_k(\rho)+B_k(\rho)|\leq \frac12|c_1^k|\rho^{\sigma_k^+}+b_k \rho^{\sigma_k^+/2},
\end{align}
 where
$$
b_k=C_2\|u\|_{H^1(B(0,1))}\bigg(\frac{\|h\|_{L^{p}(B(0,1))}}
{\omega_N(\sigma_k^+-\sigma_k^-)}\bigg)^{\!j_k}
|\mu_k|^{\big\lfloor\!\frac{N-1}4\!\big\rfloor+1}\,\frac{\big(\frac2\e\big)^{j_k-1}}
{(j_k-1)!}.
$$
From (\ref{eq:weyl}), $\mu_k\sim k^{2/(N-1)}$ and $j_k\sim
k^{1/(N-1)}$ as $k\to+\infty$. Hence we have that 
\begin{equation}\label{eq:19}
|b_k|\leq
C_3\|u\|_{H^1(B(0,1))}\exp\big(-C_4\, k^{1/(N-1)}\big),
\end{equation}
for some positive constants $C_3$ and $C_4$ depending
only on $N$, $h$, and $a$.
In view of (\ref{eq:11}) and (\ref{eq:29}), we deduce that
$$
\frac12|c_1^k|\rho^{\sigma_k^+}\leq|\varphi_k(\rho)|
+|b_k|\rho^{\sigma_k^+/2},\quad\text{for all }k\geq\bar k,
$$
and consequently
$$
\frac14\sum_{k=\bar k}^{\infty}|c_1^k|^2\int_0^1\rho^{2\sigma_k^++N-1}\,d\rho
\leq2\sum_{k=\bar k}^{\infty}\int_0^1|\varphi_k(\rho)|^2\rho^{N-1}\,d\rho
+2\sum_{k=\bar k}^{\infty} |b_k|^2\int_0^1\rho^{\sigma_k^++N-1}\,d\rho.
$$
Hence, from (\ref{eq:parseval}) and (\ref{eq:19}), we obtain that
\begin{equation}\label{eq:22}
  \sum_{k=\bar k}^{\infty}\frac{|c_1^k|^2}{N+2\sigma_k^+}
\leq 8 \|u\|^2_{L^2(B(0,1))}+8 \sum_{k=\bar k}^{\infty}\frac{|b_k|^2}{N+\sigma_k^+}
<+\infty.
\end{equation}
From \eqref{sviluppo}, 
\begin{equation}\label{eq:26}
u(\rho\,\theta)\rho^{-\sigma^+_1}=
\rho^{-\sigma^+_1}\varphi_1(\rho)\psi_1(\theta)+\sum_{k=2}^{\bar
  k-1}\rho^{\sigma_k^+-\sigma^+_1}\rho^{-\sigma^+_k}\varphi_k(\rho)\psi_k(\theta)+\sum_{k=\bar
  k}^{\infty}\rho^{-\sigma^+_1}\varphi_k(\rho)\psi_k(\theta).
\end{equation}
From (\ref{eq:20}) we deduce that
\begin{equation}\label{eq:27}
\lim_{\rho\to 0^+}\rho^{-\sigma^+_1}\varphi_1(\rho)\psi_1(\theta)=
\bigg[c_1^1+ \int_0^1
\frac{s^{-\sigma^+_1+1}}{\sigma^+_1-\sigma^-_1} h(s)\varphi_1(s)\,ds\bigg]\psi_1(\theta)
\end{equation}
and 
\begin{equation}\label{eq:28}
\lim_{\rho\to 0^+}\sum_{k=2}^{\bar
  k-1}\rho^{\sigma_k^+-\sigma^+_1}\rho^{-\sigma^+_k}\varphi_k(\rho)\psi_k(\theta)=0.
\end{equation}
From (\ref{eq:11}), (\ref{eq:29}), (\ref{eq:22}), (\ref{eq:19}), and
(\ref{eq:25}), we deduce that there exists some positive constant
$C_5$ depending only on $N$, $h$, and $a$, such that, for all
$\rho\in(0,1/2)$,
\begin{align}\label{eq:49}
  \sum_{k=\bar
    k}^{\infty}&\rho^{-\sigma^+_1}|\varphi_k(\rho)||\psi_k(\theta)|\leq
  C_1\sum_{k=\bar
    k}^{\infty}\rho^{-\sigma^+_1}\bigg(\frac32|c_1^k|\rho^{\sigma^+_k}+|b_k|
  \rho^{\sigma^+_k/2}\bigg) |\mu_k|^{\lfloor (N-1)/4\rfloor+1}\\
  \notag&\leq \frac32\,C_1\rho^{\sigma^+_2-\sigma^+_1}\bigg(\sum_{k=\bar
    k}^{\infty}\frac{|c_1^k|^2}{N+2\sigma_k^+}\bigg)^{\!\frac 12}
  \bigg(\sum_{k=\bar
    k}^{\infty}(N+2\sigma_k^+)\rho^{2(\sigma^+_k-\sigma^+_2)}|\mu_k|^{2\lfloor
    (N-1)/4\rfloor+2}\bigg)^{\!\frac 12} \\
&\notag\quad+C_1\rho^{(\sigma^+_{\bar
      k}/2)-\sigma^+_1}\sum_{k=\bar
    k}^{\infty}|b_k| |\mu_k|^{\lfloor (N-1)/4\rfloor+1}\\
  &\notag\leq C_5\|u\|_{H^1(B(0,1))}\big(\rho^{2(\sigma^+_k-\sigma^+_2)}+\rho^{(\sigma^+_{\bar
      k}/2)-\sigma^+_1}\big),
\end{align}
which implies
\begin{equation}\label{eq:30}
\lim_{\rho\to 0^+}\sum_{k=\bar
  k}^{\infty}\rho^{-\sigma^+_1}\varphi_k(\rho)\psi_k(\theta)=0.
\end{equation}
Collecting (\ref{eq:26}), (\ref{eq:27}), (\ref{eq:28}), and
(\ref{eq:30}), we finally obtain that
\begin{equation}\label{eq:31}
\lim_{\rho\to 0^+}u(\rho\,\theta)\rho^{-\sigma^+_1}= \bigg[c_1^1+
\int_0^1 \frac{s^{-\sigma^+_1+1}}{\sigma^+_1-\sigma^-_1}
h(s)\varphi_1(s)\,ds\bigg]\psi_1(\theta).
\end{equation}
We notice that, in view of (\ref{eq:42}), (\ref{eq:3}), and (\ref{eq:43}),
\begin{align}\label{eq:44}
 c_1^1&+
\int_0^1 \frac{s^{-\sigma^+_1+1}}{\sigma^+_1-\sigma^-_1}
h(s)\varphi_1(s)\,ds\\
\notag&=\varphi_1(1)-c_2^1+ \int_0^1
  \frac{s^{-\sigma^+_1+1}}{\sigma^+_1-\sigma^-_1}
  h(s)\varphi_1(s)\,ds\\
  \notag&=\int_{{\mathbb
      S}^{N-1}}u(\eta)\psi_1(\eta)\,dV(\eta) +\int_{{\mathbb
      S}^{N-1}} \bigg[ \int_0^1 \frac{s \,
    \big(s^{-\sigma^+_1}-s^{-\sigma^-_1}\big)}
  {\sigma^+_1-\sigma^-_1}\,h(s)u(s\,\eta)\,ds\bigg]\psi_1(\eta)
  \,dV(\eta).
\end{align}
The limit in (\ref{eq:40}) follows now from (\ref{eq:31}) and
(\ref{eq:44}) in the case $r=1$ and by a change of variable inside the
integral for $r\neq1$.  Moreover estimates (\ref{eq:46}),
(\ref{eq:23}), (\ref{eq:24}), (\ref{eq:47}), (\ref{eq:48}), the
definition of $j_k$, and (\ref{eq:49}), imply that, for some $C_6>0$
depending only on $N$, $R$, $h$, and $a$,
\begin{align}\label{eq:51}
  u(\rho\,\theta)\rho^{-\sigma^+_1}&=\bigg(c_1^1+
  \int_\rho^1\frac{s^{-\sigma^+_k+1}}{\sigma^+_k-\sigma^-_k}
  h(s)\varphi_1(s)\bigg)\psi_1(\theta)+\rho^{-\sigma^+_1}\bigg[
  B_1(\rho)\psi_1(\theta)+\sum_{k>1}^{\infty}
\varphi_k(\rho)\psi_k(\theta)\bigg]\\
&\notag\leq C_6\|u\|_{H^1(B(0,R))},
\end{align}
for all $0<\rho<1/2$. On the other hand, standard elliptic estimates
in $B(0,1)\setminus B(0,1/2)$ yield, for some $C_7>0$ depending only on
$N$, $R$, $h$, and $a$,
\begin{align}\label{eq:52}
  u(\rho\,\theta)\rho^{-\sigma^+_1}\leq C_7\|u\|_{H^1(B(0,R))},
\quad\text{for all }\frac12\leq\rho<1.
\end{align}
Estimate (\ref{eq:50}) follows from (\ref{eq:51}) and (\ref{eq:52}).

From (\ref{eq:31}) and the positivity of $u$, it follows easily that
$\Big[c_1^1+ \int_0^1 \frac{s^{-\sigma^+_1+1}}{\sigma^+_1-\sigma^-_1}
h(s)\varphi_1(s)\,ds\Big]\geq 0$. Since 
$$0
<\min_{{\mathbb
    S}^{N-1}}\psi_1\leq \psi_1(\theta)\leq \max_{{\mathbb
    S}^{N-1}}\psi_1,\quad\text{for all }\theta\in{\mathbb
    S}^{N-1},
$$
and, by standard regularity theory, $u\in
C^0\big(\overline{B(0,1)}\setminus B(0,s)\big)$ for any $s\in(0,1)$,
the proof of Proposition \ref{p:1} will be complete if we show that
\begin{equation}\label{eq:32}
c_1^1+
\int_0^1 \frac{s^{-\sigma^+_1+1}}{\sigma^+_1-\sigma^-_1}
h(s)\varphi_1(s)\,ds>0.
\end{equation}
In order to obtain (\ref{eq:32}), we need to prove the following

\medskip\noindent {\bf Claim 3:} if $k\in \N\setminus\{0\}$ and 
$$
c_1^k+\int_0^1
\frac{s^{-\sigma^+_k+1}}{\sigma^+_k-\sigma^-_k} h(s)\varphi_k(s)\,ds=0.
$$
then $\varphi_k(\rho)=0$ for all $\rho\in(0,1)$. Indeed, if $c_1^k+ \int_0^1
\frac{s^{-\sigma^+_k+1}}{\sigma^+_k-\sigma^-_k} h(s)\varphi_k(s)\,ds=0$, then 
\begin{equation}\label{eq:33}\varphi_k(\rho)=
-\rho^{\sigma^+_k}
\int_0^\rho\frac{s^{-\sigma^+_k+1}}{\sigma^+_k-\sigma^-_k}
h(s)\varphi_k(s)\,ds+
\rho^{\sigma^-_k}
\int_0^\rho\frac{s^{-\sigma^-_k+1}}{\sigma^+_k-\sigma^-_k}
h(s)\varphi_k(s)\,ds.
\end{equation}
From claim 2 and (\ref{eq:20}), we know that there exists a constant $\ell_k$  
depending on $k$ (and on $a$, $R$, $h$, $u$, $N$) such  that
$$
|\varphi_k(\rho)|\leq\ell_k\rho^{\sigma^+_k}.
$$
Using the above estimate in (\ref{eq:33}), we can improve such an estimate as
$$
|\varphi_k(\rho)|\leq\ell_k\rho^{\sigma^+_k}\frac{2\|h\|_{L^{p}(B(0,1))}}
{\omega_N(\sigma_k^+-\sigma_k^-)}\frac{\rho^{\e}}{\e}.
$$
Using the above estimate in (\ref{eq:33}), we can obtain the following further improvement
$$
|\varphi_k(\rho)|\leq \ell_k
\rho^{\sigma^+_k}\bigg(\frac{2\|h\|_{L^{p}(B(0,1))}}
{\omega_N(\sigma_k^+-\sigma_k^-)\e}\bigg)^2\frac{(\rho^{\e})^2}{(2\cdot1)^{\frac{p-1}{p}}}.
$$
Arguing by induction, we can easily prove that, for all $j\in\N$,
$$
|\varphi_k(\rho)|\leq \ell_k
\rho^{\sigma^+_k}\bigg(\frac{2\|h\|_{L^{p}(B(0,1))}}
{\omega_N(\sigma_k^+-\sigma_k^-)\e}\bigg)^j\frac{(\rho^{\e})^j}{(j!)^{\frac{p-1}{p}}},
$$
and letting $j\to+\infty$, we deduce that $\varphi_k(\rho)=0$ for all
$\rho\in(0,1)$. Claim 3 is thereby proved.

\medskip\noindent We are now in position to prove (\ref{eq:32}).
Arguing by contradiction, let us assume that 
\begin{equation}\label{eq:34}
c_1^1+
\int_0^1 \frac{s^{-\sigma^+_1+1}}{\sigma^+_1-\sigma^-_1}
h(s)\varphi_1(s)\,ds=0
\end{equation}
and let $k_0>1$ be the smallest index for which 
$$
c_1^{k_0}+ \int_0^1
\frac{s^{-\sigma^+_{k_0}+1}}{\sigma^+_{k_0}-\sigma^-_{k_0}} h(s)\varphi_{k_0}(s)\,ds\not=0.
$$
Such a $k_0$ exists in view of claim 3; indeed if $c_1^k+ \int_0^1
\frac{s^{-\sigma^+_k+1}}{\sigma^+_k-\sigma^-_k}
h(s)\varphi_k(s)\,ds=0$ for all $k$, then $\varphi_k\equiv0$ for all
$k$ and $u$ would be identically zero, thus giving rise to a
contradiction.  Moreover, from (\ref{eq:34}), we have that $k_0>1$,
and, by claim 3, $\varphi_k\equiv0$ in $(0,1)$ for all $1\leq k\leq
k_0-1$. Repeating the same arguments we used above to prove
(\ref{eq:31}), it is now possible to show that
\begin{equation}\label{eq:55}
\lim_{\rho\to 0^+}u(\rho\,\theta)\rho^{-\sigma^+_{k_0}}= \sum_{k=k_0}^{k_0+m_{k_0}-1}
\bigg[c_1^{k}+
\int_0^1 \frac{s^{-\sigma^+_{k}+1}}{\sigma^+_{k}-\sigma^-_{k}}
h(s)\varphi_{k}(s)\,ds\bigg]\psi_{k}(\theta),
\end{equation}
where $m_{k_0}$ is the geometric multiplicity of the eigenvalue
$\mu_{k_0}$. We notice that the sum at the right hand side is a nontrivial 
function in $L^2({\mathbb S}^{N-1})$ which, being $k_0>1$, is orthogonal to the 
first positive eigenfunction $\psi_1$. Hence the  right hand side of \eqref{eq:55}
 changes sign in
${\mathbb S}^{N-1}$. Therefore  the limit in \eqref{eq:55} 
implies that $u$ changes sign in
a neighborhood of $0$, which is in contradiction with the positivity
assumption on $u$. Condition (\ref{eq:32}) follows and the proof of
Proposition \ref{p:1} is now complete.
  \end{pf}

\begin{remark}\label{r:changingsign}
If we let the assumption of positivity of $u$ drop, following the proof of 
Proposition~\ref{p:1}, we can still prove a {\it Cauchy's integral type formula}
for $u$. More precisely, if 
$u\in H^1(B(0,R))$ is a weak solution to
\eqref{eq:1} in $B(0,R)$ which changes sign in any neighborhood of $0$, 
with  a radial potential $h\in  L^{\infty}_{\rm loc}(0,R)\cap L^p(0,R)$ for
some $p>N/2$, then, following the notation introduced in~\eqref{eq:56} and
 letting  $k_0>1$ be the smallest index for which 
\begin{multline*}
\int_{{\mathbb
      S}^{N-1}} \bigg(r^{-\sigma_{k_0}^+} u(r\eta)+\int_0^r 
{{\frac{s^{1-\sigma_{k_0}^+}}
  {\sigma_{k_0}^+-\sigma_{k_0}^-}}}\,h(s)u(s\,\eta)\,ds
\\
-r^{\sigma_{k_0}^--\sigma_{k_0}^+}\int_0^r 
{{\frac{s^{1-\sigma_{k_0}^-}}
  {\sigma_{k_0}^+-\sigma_{k_0}^-}}}\,h(s)u(s\,\eta)\,ds
\bigg)\psi_{k_0}(\eta)
  \,dV(\eta)\not=0,
\end{multline*}
for any
$\theta\in{{\mathbb S}^{N-1}}$ and $r\in(0,R)$ there holds
\begin{align*}
\lim_{\rho\to0^+}u(\rho\,\theta)
\rho^{-\sigma_{k_0}^+}&=\sum_{\{k:\,\mu_k=\mu_{k_0}\}}
\bigg[\int_{{\mathbb
      S}^{N-1}} \bigg(r^{-\sigma_k^+} u(r\eta)+\int_0^r 
{{\frac{s^{1-\sigma_k^+}}
  {\sigma_k^+-\sigma_k^-}}}\,h(s)u(s\,\eta)\,ds
\\
&\hskip3.5cm-r^{\sigma_k^--\sigma_k^+}\int_0^r 
{{\frac{s^{1-\sigma_k^-}}
  {\sigma_k^+-\sigma_k^-}}}\,h(s)u(s\,\eta)\,ds
\bigg)\psi_k(\eta)
  \,dV(\eta)\bigg]\psi_{k}(\theta).
\end{align*}
\end{remark}

Without the assumption of radial symmetry of the potential, it is
still possible to evaluate the exact behavior near the singularity of
the first Fourier coefficient $\varphi_1$ (see \eqref{sviluppo} and (\ref{eq:3})).

\begin{Lemma}\label{l:limite}
Let $a\in L^{\infty}({\mathbb S}^{N-1})$ be such that
  $\Lambda_N(a)<1$, $R>0$, and $u\in H^1(B(0,R))$, $u\geq 0$ a.e. in
  $B(0,R)$, $u\not\equiv 0$, be a weak $H^1$-solution to
\begin{equation*}
  -\Delta u(x)=\bigg[\frac{a(x/|x|)}
  {|x|^2}+q(x)\bigg]\,u(x)\quad\text{ in }B(0,R),
\end{equation*} 
where $q\in L^{\infty}_{\rm loc}\big(B(0,R)\setminus \{0\}\big)\cap
L^p(B(0,R))$ for some $p>\frac N2$. Then, for any $0<r<R$, 
\begin{gather}\label{eq:limit}
\lim_{\rho\to0^+}\rho^{-\sigma}\int_{{\mathbb S}^{N-1}}u(\rho\,\theta)
\psi_1(\theta)\,dV(\theta)=\int_{{\mathbb
S}^{N-1}}\bigg(r^{-\sigma}u(r\,\theta)+\int_0^r 
{{\frac{s^{1-\sigma}}
{2\sigma+N-2}}}\,q(s\,\theta)u(s\,\theta)\,ds
\\
\notag-r^{-2\sigma-N+2}\int_0^r 
{{\frac{s^{N-1+\sigma}}
  {2\sigma+N-2}}}\,q(s\,\theta)u(s\,\theta)\,ds
\bigg)\psi_1(\theta)
  \,dV(\theta).
\end{gather}

\end{Lemma}

\begin{pf}
  The proof follows the lines of the first part of the proof of
  Proposition \ref{p:1}. By scaling, it is sufficient to prove 
(\ref{eq:limit}) for $r=1$. Let
\begin{equation*}
u(x)=u(\rho\,\theta)=\sum_{k=1}^\infty\varphi_k(\rho)\psi_k(\theta)
\quad\mbox{ and }\quad 
q(x)u(x)=q(\rho\,\theta)u(\rho\,\theta)=
\sum_{k=1}^\infty\zeta_k(\rho)\psi_k(\theta)
\end{equation*}
where $\rho=|x|\in(0,1]$, $\theta=x/|x|\in{{\mathbb S}^{N-1}}$, 
\begin{equation*}
\varphi_k(\rho)=\int_{{\mathbb S}^{N-1}}u(\rho\,\theta)\psi_k(\theta)\,dV(\theta)\,,
\quad
\zeta_k(\rho)=\int_{{\mathbb S}^{N-1}}
q(\rho\,\theta)u(\rho\,\theta)\psi_k(\theta)\,dV(\theta),
\end{equation*}
and $\psi_k$ is an $L^2$-normalized eigenfunction of the operator $-\D_{\mathbb
S^{N-1}}-a(\theta)$ on the sphere associated to
the $k$-th eigenvalue $\mu_k$, i.e. satisfying \eqref{eq:2rad}.
The first Fourier coefficient $\varphi_1$ solves 
\begin{equation*}
\varphi_1''(\rho)+\frac{N-1}{\rho}\,\varphi_1^\prime(\rho)-
\frac{\mu_1}{\rho^2}\,\varphi_1(\rho)=\zeta_1(\rho)\quad\text{in }(0,1).
\end{equation*}
A direct calculation shows that, for some $c_1^1,c_2^1\in\R$,
\begin{equation*}
\varphi_1(\rho)=\rho^{\sigma^+_1}
\bigg(c_1^1+\int_\rho^1\frac{s^{-\sigma^+_1+1}}{\sigma^+_1-\sigma^-_1}
\zeta_1(s)\,ds\bigg)+\rho^{\sigma^-_1}
\bigg(c_2^1+\int_\rho^1\frac{s^{-\sigma^-_1+1}}{\sigma^-_1-\sigma^+_1}
\zeta_1(s)\,ds\bigg),
\end{equation*}
where $\sigma^+_1=-\frac{N-2}{2}+\sqrt{\big(\frac{N-2}
  {2}\big)^2+\mu_1}$ and
$\sigma^-_1=-\frac{N-2}{2}-\sqrt{\big(\frac{N-2}{2}\big)^2+\mu_1}.$
From Theorem \ref{t:bk} and standard regularity theory, we deduce
that $u(x)\leq {\rm const\,}|x|^{\sigma_1^+}$ in $B(0,1)$, hence, 
 H\"older's inequality yields
\begin{align}\label{eq:limit1}
\int_0^1\big|s^{-\sigma^+_1+1}
\zeta_1(s)\big|\,ds&\leq
\int_0^1s^{-\sigma^+_1+1}
\bigg(\int_{{\mathbb
S}^{N-1}}|q(s\,\theta)|u(s\,\theta)\psi_1(\theta)\,dV(\theta)\bigg)\,ds\\[5pt]
\notag&\leq {\rm const\,}\,\int_0^{1}s^{1-\frac{N-1}{p}}
\bigg(\int_{{\mathbb
S}^{N-1}}s^{\frac{N-1}{p}}|q(s\,\theta)|\,dV(\theta)\bigg)\,ds\\[5pt]
\notag&\leq {\rm const\,}\|q\|_{L^{p}(B(0,1))}
\bigg(\int_0^1
s^{(1-\frac{N-1}{p})
\frac{p}{p-1}}\,ds\bigg)^{1-\frac1p}<\infty.
\end{align}
In a similar way, we obtain that 
\begin{equation}\label{eq:7limit}
s\mapsto \frac{s^{-\sigma^-_1+1}}{\sigma^-_1-\sigma^+_1}
\,\zeta_1(s)\in L^1(0,1).
\end{equation}
Since $u\in L^{2^*}(B(0,1))$, $\sigma^-_1<\sigma^+_1$, 
from (\ref{eq:limit1}),
(\ref{eq:7limit}), and the fact that $\rho^{\sigma^-_1}\psi_1(\theta)\not
\in L^{2^*}(B(0,1))$, we conclude that there must be
\begin{equation*}
c_2^1=-\int_0^1\frac{s^{-\sigma^-_1+1}}{\sigma^-_1-\sigma^+_1}
\,\zeta_1(s)\,ds,
\end{equation*}
hence 
\begin{align}\label{eq:8limit}
\rho^{-\sigma^+_1}\varphi_1(\rho)=
c_1^1+\int_\rho^1\frac{s^{-\sigma^+_1+1}}{\sigma^+_1-\sigma^-_1}
\,\zeta_1(s)\,ds+
\rho^{\sigma^-_1-\sigma^+_1}\int_0^\rho\frac{s^{-\sigma^-_1+1}}{\sigma^+_1-\sigma^-_1}
\,\zeta_1(s)\,ds\quad\text{for any }\rho\in(0,1).
\end{align}
Notice that, by standard regularity theory, 
$\varphi_1$ is continuous at $\rho=1$, thus, letting $\rho\to 1^-$ in 
(\ref{eq:8limit}), we obtain  
\begin{equation}\label{eq:59limit}
c_1^1=\varphi_1(1)-
\int_0^1\frac{s^{-\sigma^-_1+1}}{\sigma^+_1-\sigma^-_1}
\,\zeta_1(s)\,ds.
\end{equation}
Arguing as in \eqref{eq:limit1}, we obtain that
\begin{equation}\label{eq:61limit}
\rho^{\sigma^-_1-\sigma^+_1}\int_0^\rho\Big|\frac{s^{-\sigma^-_1+1}}{\sigma^-_1-\sigma^+_1}
\zeta_1(s)\Big|\,ds\leq {\rm const\,}\rho^{2-\frac{N}{p}}.
\end{equation}
Since $p>\frac{N}{2}$, \eqref{eq:8limit}, \eqref{eq:59limit}, and \eqref{eq:61limit}
imply that
\begin{align*}
\lim_{\rho\to0^+}\rho^{-\sigma^+_1}\varphi_1(\rho)=
\varphi_1(1)-\int_0^1
\frac{s^{-\sigma^-_1+1}}{\sigma^+_1-\sigma^-_1}
\,\zeta_1(s)\,ds+\int_0^1\frac{s^{-\sigma^+_1+1}}{\sigma^+_1-\sigma^-_1}
\,\zeta_1(s)\,ds,
\end{align*}
and \eqref{eq:limit} for $r=1$ follows. The result in the case  $r\neq1$
 can be easily obtained  just by scaling.~\end{pf}

\noindent In order to extend the result of Proposition \ref{p:1} to
the case in which the potential is a non radial perturbation of the
dipole-type singular term, we will construct a subsolution and a
supersolution  which solve equations of
type (\ref{eq:1}) and the behavior of which is consequently known in
view of Proposition \ref{p:1}.

\begin{Lemma}\label{l:ulsol}
Let  $a\in L^{\infty}({\mathbb S}^{N-1})$ be such that
  $\Lambda_N(a)<1$, $C\in\R$, and $\e>0$. Then, for all 
\begin{equation}\label{eq:2}
0<r<
\begin{cases}
  \Big[\frac{(N-2)^2}{4C^+}(1-\Lambda_N(a))\Big]^{\!1/\e},&\text{if }C>0,\\[5pt]
+\infty,&\text{if }C\leq0,
\end{cases}
\end{equation}
and for all $\gamma\in H^{1/2}(\partial B(0,r))$, 
$\gamma\geq0$, $\gamma\not\equiv0$,  the Dirichlet boundary value problem  
\begin{equation}\label{eq:uleq}
\left\{\begin{array}{ll}
-\Delta u(x)=\bigg[\dfrac{a(x/|x|)}{|x|^2}
+C|x|^{-2+\e}\bigg]\,u(x),& \text{ in }B(0,r),\\[10pt]
u\big|_{\partial B(0,r)}=\gamma, & \text{ on }\partial B(0,r),\\
\end{array}\right.
\end{equation}
admits a unique weak solution $u\in H^1(B(0,r))$. Moreover $u$ is
continuous  and strictly positive in $B(0,r)\setminus\{0\}$, and there
exists a positive constant $C'$ depending on $a$, $C$, $\e$, $N$, and
$r$, such that
 \begin{align}\label{eq:53}
\|u\|_{H^1(B(0,r))}\leq C'\|\gamma\|_{H^{1/2}(\partial B(0,r))}.
\end{align}
In addition, if $\gamma\in W^{2-1/k,k}(\partial B(0,r))$ for some $k>N/2$, then 
$u\in C^0(\overline{B(0,r)}\setminus\{0\})$.
\end{Lemma}
\begin{pf}
  For a fixed $r$ satisfying (\ref{eq:2}) and $\gamma\in H^{1/2}(\partial B(0,r))$, 
$\gamma\geq0$, $\gamma\not\equiv0$, let $\tilde v$ be the
unique $H^1(B(0,r))$-weak solution to the problem
$$
\begin{cases}
\tilde v\in H^1(B(0,r)),\\
-\Delta \tilde v=0,&\text{in }B(0,r),\\
v=\gamma,&\text{on }\partial B(0,r).
\end{cases}
$$
By classical trace embedding theorems, it follows that 
\begin{equation}\label{eq:41}
\|\tilde v\|_{H^1(B(0,r))}
\leq {\rm const}(N,r)\|\gamma\|_{H^{1/2}(\partial B(0,r))},
\end{equation}
for some positive constant ${\rm const}(N,r)$ depending only on 
$N$ and $r$.
Let us define the quadratic form $a:H^1_0(B(0,r))\times H^1_0(B(0,r))\to\R$ as 
\begin{align*}
  {\mathcal Q}(w,u):=\int_{B(0,r)}\bigg[ \n w(x)\cdot\nabla
  w(x)-\frac{1}{|x|^2}\big(a(x/|x|)+C|x|^{\e}\big)w(x)u(x)\bigg]\,dx,
\end{align*}
and $\Phi\in H^{-1}(B(0,r))$ as
$${}_{H^{-1}}\big\langle
\Phi,u\big\rangle_{H^1_0}:=\int_{B(0,r)}\bigg(\frac{a(x/|x|)}{|x|^2}
+\frac{C}{|x|^{2-\e}}\bigg)\tilde v(x)u(x)\,dx.
$$
By Hardy's inequality, it is easy to verify that 
\begin{equation}\label{eq:35}
{\mathcal Q}(u,u)\geq\bigg[1-\Lambda_N(a)-\frac{4C^+ r^{\e}}{(N-2)^2}
\bigg]
\int_{B(0,r)}|\n u(x)|^2\,dx.
\end{equation}
Since (\ref{eq:2}) implies that
$\big[1-\Lambda_N(a)-\frac{4C^+ r^{\e}}{(N-2)^2}
\big]>0$, we conclude that
the bilinear bounded form ${\mathcal Q}$ is
coercive. Furthermore, the function $x\mapsto a(x/|x|)|x|^{-2}+C|x|^{-2+\e}$
 belongs to $L^{\frac{2N}{N+2}}(B(0,r))$,
hence $\Phi$ is a bounded linear functional on $H^1_0(B(0,r))$. From
the Lax-Milgram lemma we deduce that there exists a unique $w\in
H^1_0(B(0,r))$ such that ${\mathcal Q}(w,u)={}_{H^{-1}}\big\langle
\Phi,u\big\rangle_{H^1_0}$ for all $u\in H^1_0(B(0,r))$. 
In particular $w$ weakly solves
\begin{equation}\label{eq:77}
\begin{cases}
-\Delta w(x)-\dfrac{1}{|x|^2}\Big[a(x/|x|)+C|x|^{\e}
\Big]w(x)=\bigg[\dfrac{a(x/|x|)}{|x|^2}+
\dfrac{C}{|x|^{2-\e}}\bigg]\tilde v(x),&\text{in }B(0,r),\\
w=0,&\text{on }\partial B(0,r).
\end{cases}
\end{equation}
Testing the above equation with $w$ and using (\ref{eq:35}),
Poincar\'e's and H\"older's inequalities and \eqref{eq:41}, we obtain that
\begin{align*}
\|w\|_{H^1(B(0,r))}\leq c(a,C,\e,N,r)\|\tilde v\|_{H^1(B(0,r))}
\leq c'(a,C,\e,N,r)\|\gamma\|_{H^{1/2}(\partial B(0,r))},
\end{align*}
for some positive constants $c(a,C,\e,N,r)$ and $c'(a,C,\e,N,r)$
depending on $a$, $C$, $\e$, $N$, and $r$.  It is now easy to verify
that $u:=w+\tilde v\in H^1(B(0,r))$ satisfies (\ref{eq:53}) and is the
unique weak solution to (\ref{eq:uleq}). Moreover, testing
(\ref{eq:uleq}) with $-u^-:=-\max\{-u,0\}$ and using (\ref{eq:35}), we
obtain that
$$
0={\mathcal Q}(u^-,u^-)\geq\bigg[1-\Lambda_N(a)-\frac{4C^+ r^{\e}}{(N-2)^2}
\bigg]
\int_{B(0,r)}|\n u^-(x)|^2\,dx,
$$
which, in view of (\ref{eq:2}), implies that $u^-=0$ a.e. in $B(0,r)$,
i.e. $u\geq 0$ a.e in $B(0,r)$. The Strong Maximum Principle allows us
to conclude that $u>0$ in $B(0,r)\setminus\{0\}$, while standard
regularity theory for elliptic equations ensures interior continuity
of $u$ outside the origin.

If, in addition, we assume that $\gamma\in W^{2-1/k,k}(\partial
B(0,r))$ for some $k>N/2$, then $\tilde v\in W^{2,k}(B(0,r)$, and
hence $\tilde v\in C^{0,\a}(\overline{B(0,r)})$, hence, from elliptic
regularity theory applied to (\ref{eq:77}) outside $0$, we obtain that $u\in
C^0(\overline{B(0,r)}\setminus\{0\})$.
\end{pf}

\begin{pfn}{Theorem \ref{t:1}}
  Let $R>0$ such that $\overline{B(0,R)}\subset\Omega$.  Since
  $q(x)=O(|x|^{-(2-\e)})$ as $|x|\to0$ for some $\e>0$ and $q\in
  L^{\infty}_{\rm loc}\big(\Omega\setminus \{0\}\big)$, there exists a
  positive constant $\tilde C$ such that $-\tilde C |x|^{-(2-\e)}\leq
  q(x)\leq \tilde C |x|^{-(2-\e)}$ for a.e. $x\in B(0,R)$. Let us fix
  $\bar r=\bar r(R,N,q,a,\e)$, such that $0<\bar r<\min\Big\{R,
  \big[\frac{(N-2)^2}{4\tilde C}(1-\Lambda_N(a))\big]^{\!1/\e}
  \Big\}$.  We notice that the Maximum Principle implies that $u>0$ in
  $\overline{B(0,R)}\setminus\{0\}$, whereas standard elliptic
  regularity theory yields $u\in W^{2,k}(\overline{B(0,R)} \setminus
  B(0,s))$ for all $s\in(0,R)$ and some $k>N/2$, and, consequently,
 $u$ is continuous in $\overline{B(0,R)}\setminus\{0\}$.  Hence the function
  $\gamma_r:=u\big|_{\partial B(0,r)}$ belongs to $W^{2-1/k,k}(\partial
B(0,r))$ for some $k>N/2$ and is continuous and strictly positive on $\partial
  B(0,r)$ for all $0<r\leq \bar r$.  From Lemma \ref{l:ulsol} we
  deduce that, for any $0<r\leq \bar r$, there exist $\underbar{\it
    u}_r\in H^1(B(0,r))$ and $\bar u^r\in H^1(B(0,r))$ continuous and
  strictly positive in $\overline{B(0,r)}\setminus\{0\}$, weakly satisfying
\begin{equation*}
\begin{cases}
-\Delta \underbar{\it u}_r(x)=\bigg[\dfrac{a(x/|x|)}{|x|^2}
-\tilde C|x|^{-2+\e}\bigg]\,\underbar{\it u}_r(x),& \text{ in }B(0,r),\\[10pt]
\underbar{\it u}_r\big|_{\partial B(0,r)}=\gamma_r, & \text{ on }\partial B(0,r),
\end{cases}
\end{equation*}
and 
\begin{equation*}
\begin{cases}
-\Delta \bar{u}_r(x)=\bigg[\dfrac{a(x/|x|)}{|x|^2}
+\tilde C|x|^{-2+\e}\bigg]\,\bar{u}_r(x),& \text{ in }B(0,r),\\[10pt]
\bar{u}_r\big|_{\partial B(0,r)}=\gamma_r, & \text{ on }\partial B(0,r).
\end{cases}
\end{equation*}
From Proposition \ref{p:1}, there exist two constants
$A_2>A_1>0$ (depending on $\bar r$, $N$, $q$, $a$, $\e$, and $u$) such that
\begin{equation}\label{eq:38}
A_1|x|^{\sigma}\leq 
\underbar{\it u}_{\bar r}(x)\quad\text{and}\quad
\bar{u}_{\bar r}(x)\leq
A_2|x|^{\sigma},\quad\text{ for all }x\in B(0,\bar r/2)\setminus\{0\}.
\end{equation}
Furthermore, for all $0<r\leq\bar r$, $u-\underbar{\it u}_r$ satisfies
\begin{equation}\label{eq:36}
\begin{cases}
  -\Delta (u-\underbar{\it
    u}_r)(x)-\bigg[\dfrac{a(x/|x|)}{|x|^2}
  -\tilde C|x|^{-2+\e}\bigg]\,(u-\underbar{\it u}_r)(x)\geq0,
& \text{ in }B(0,r),\\[10pt]
  (u-\underbar{\it u}_r)\big|_{\partial B(0,r)}=0, & \text{ on }\partial
  B(0,r),
\end{cases}
\end{equation}
while $u-\bar{u}_r$ satisfies
\begin{equation}\label{eq:37}
\begin{cases}
  -\Delta (u-\bar{u}_r)(x)-\bigg[\dfrac{a(x/|x|)}{|x|^2}
  +\tilde C|x|^{-2+\e}\bigg]\,(u-\bar{u}_r)(x)\leq0,
& \text{ in }B(0,r),\\[10pt]
  (u-\bar{u}_r)\big|_{\partial B(0,r)}=0, & \text{ on }\partial
  B(0,r).
\end{cases}
\end{equation}
Testing (\ref{eq:36}), respectively (\ref{eq:37}), with
$-(u-\underbar{\it u}_r)^-$, respectively $(u-\bar{u}_r)^+$, and using
(\ref{eq:35}), we obtain that, for any $0<r\leq\bar r$,
\begin{equation}\label{eq:39}
\underbar{\it u}_r(x)\leq u(x)\leq \bar{u}_r(x),\quad\text{ for all }x\in
B(0,r)\setminus\{0\}.
\end{equation}
In particular, from (\ref{eq:38}) and (\ref{eq:39}), we deduce that 
\begin{equation}\label{eq:69}
A_1|x|^{\sigma}\leq 
u(x)\leq
A_2|x|^{\sigma},\quad\text{ for all }x\in B(0,\bar r/2)\setminus\{0\}.
\end{equation}
Estimate \eqref{eq:69} and the continuity of $u$ outside the origin
imply that there exists a positive constant $C$ (depending on $q$, $R$, $\Omega$,
$a$, $\e$, and $u$) such that
\begin{align}\label{eq:74}
\frac1C|x|^{\sigma}\leq u(x)\leq
C|x|^{\sigma}\quad\text{ for all }x\in \overline{B(0,R)}\setminus\{0\}.
\end{align}
Let us now fix $\delta=\delta(N,a,\e)>0$ such that
\begin{equation*}
\delta<\min\Big\{\e, {\textstyle{\sqrt{\big(\frac{N-2}
  {2}\big)^2+\mu_1}}}\Big\}
\end{equation*}
and set 
$$
\hat r=\min\Big\{\bar r,\Big[{\textstyle{\frac{\delta}{\tilde C}}}\Big(2
{\textstyle{\sqrt{\big(\frac{N-2} {2}\big)^2+\mu_1}}}-\delta\Big)\Big]^{1/\e},
\Big({\textstyle{\frac{\min_{{\mathbb S}^{N-1}}\psi_1}{C}}}\Big)^{1/\delta}\Big\},
$$
with $C$ given in \eqref{eq:74}. The function $\hat u$ defined as 
\begin{equation}\label{eq:73}
\hat u(x)=|x|^{{\sigma}-\delta}\psi_1\big(x/|x|\big)
\end{equation}
belongs to $H^1(B(0,R))$ and, for all $0<r\leq\hat r$, satisfies
$$
\begin{cases}
  -\Delta \hat{u}(x)-\dfrac{a(x/|x|)}{|x|^2}\, \hat u(x)= \delta
  \Big(2 {\textstyle{\sqrt{\big(\frac{N-2}
        {2}\big)^2+\mu_1}}}-\delta\Big) |x|^{-2}\,\hat u(x)\geq \tilde
  C\,|x|^{-2+\e}\,
  \hat u(x),& \text{ in }B(0,r),\\[10pt]
  \hat{u}(x)\big|_{\partial
    B(0,r)}=r^{\sigma-\delta}\psi_1(x/r)\geq \bar u_r(x), & \text{
    on }\partial B(0,r).
\end{cases}
$$
Hence, for all $0<r\leq\hat r$, $\hat u-\bar{u}_r$ satisfies
\begin{equation*}
\begin{cases}
  -\Delta (\hat u-\bar{u}_r)(x)-\bigg[\dfrac{a(x/|x|)}{|x|^2}\,x\cdot{\mathbf d}
  +\tilde C|x|^{-2+\e}\bigg]\,(\hat u-\bar{u}_r)(x)\geq0,
& \text{ in }B(0,r),\\[10pt]
  (\hat u-\bar{u}_r)\big|_{\partial B(0,r)}\geq0, & \text{ on }\partial
  B(0,r).
\end{cases}
\end{equation*}
Testing the above equation  with
$-(\hat u-\bar u_r)^-$ and using
(\ref{eq:35}), we obtain that, for any $0<r\leq\hat r$,
\begin{equation}\label{eq:75}
 \bar{u}_r(x)\leq \hat u(x),\quad\text{ for all }x\in
B(0,r)\setminus\{0\}.
\end{equation}
From Proposition \ref{p:1}, for any $0<r\leq \hat r$, the functions
$$
x\mapsto \frac{\underbar{\it u}_r(x)}{|x|^{\sigma}\psi_1(x/|x|)}\quad
\text{and}\quad
x\mapsto \frac{\bar u_r(x)}{|x|^{\sigma}\psi_1(x/|x|)}
$$
have limits as $|x|\to 0$, which, accordingly with
(\ref{eq:40}--\ref{eq:58}) and taking into account the continuity of functions 
$\underbar{\it
      u}_r$ and $\bar
      u_r$ up to the boundary $|x|=r$, can be computed as
\begin{align*}
  \underbar{\it L}_r:&=\lim_{|x|\to 0}\frac{\underbar{\it
      u}_r(x)}{|x|^{\sigma}\psi_1(x/|x|)}\\
&=\int_{{\mathbb
      S}^{N-1}} \bigg( r^{-\sigma} u(r\eta)
-\tilde C\int_0^{r} 
{{\frac{s^{1-\sigma}}
  {2\sigma+N-2}}}\,s^{-2+\e}\underbar{\it u}_r(s\,\eta)\,ds
\\
&\qquad+\tilde C r^{-2\sigma-N+2}\int_0^{r} 
{{\frac{s^{N-1+\sigma}}
  {2\sigma+N-2}}}\,s^{-2+\e}\underbar{\it u}_r(s\,\eta)\,ds
\bigg)\psi_1(\eta)
  \,dV(\eta),
\end{align*}
and 
\begin{align*}
  \bar L_r:&=\lim_{|x|\to 0}\frac{\bar
      u_r(x)}{|x|^{\sigma}\psi_1(x/|x|)}\\
&=\int_{{\mathbb
      S}^{N-1}} \bigg(r^{-\sigma} u(r\eta)+\tilde C\int_0^r 
{{\frac{s^{1-\sigma}}
  {2\sigma+N-2}}}\,s^{-2+\e}\bar{u}_r(s\,\eta)\,ds
\\
&\qquad-\tilde C r^{-2\sigma-N+2}\int_0^r 
{{\frac{s^{N-1+\sigma}}
  {2\sigma+N-2}}}\,s^{-2+\e}\bar{u}_r(s\,\eta)\,ds
\bigg)\psi_1(\eta)
  \,dV(\eta).
\end{align*}
From \eqref{eq:39} and \eqref{eq:74}, it follows that
$$
\underbar{\it L}_r=r^{-\sigma}\int_{{\mathbb
      S}^{N-1}} u(r\eta)\psi_1(\eta)
  \,dV(\eta)+o(1)\quad \text{as }r\to 0.
$$
From \eqref{eq:75}, \eqref{eq:73}, and the choice of $\delta$, we obtain that
$$
  \bar L_r=r^{-\sigma}\int_{{\mathbb
      S}^{N-1}} u(r\eta)\psi_1(\eta)
  \,dV(\eta)+o(1)\quad \text{as }r\to 0.
$$
Hence, from Lemma \ref{l:limite}, we conclude that, for any $R$ such that 
$\overline{B(0,R)}\subset\Omega$, 
\begin{align}\label{eq:76}
\lim_{r\to 0}\underbar{\it L}_r=
\lim_{r\to 0}  \bar L_r=&
\int_{{\mathbb
S}^{N-1}}\bigg(R^{-\sigma}u(R\,\eta)+\int_0^R 
{{\frac{s^{1-\sigma}}
{2\sigma+N-2}}}\,q(s\,\eta)u(s\,\eta)\,ds
\\
&\notag-R^{-2\sigma-N+2}\int_0^R 
{{\frac{s^{N-1+\sigma}}
  {2\sigma+N-2}}}\,q(s\,\eta)u(s\,\eta)\,ds
\bigg)\psi_1(\eta)
  \,dV(\eta).
\end{align}
In view of \eqref{eq:39}, there holds that, for any $0<r\leq\bar r$,
{\begin{align*}
\underline{ L}_r&=\lim_{|x|\to 0}\frac{\underline{u}_r(x)}{|x|^{a_{\mu_1}}\psi_1(x/|x|)}\leq
\liminf_{|x|\to 0}\frac{u(x)}{|x|^{a_{\mu_1}}\psi_1(x/|x|)}\\
& \leq \limsup_{|x|\to 0}\frac{u(x)}{|x|^{a_{\mu_1}}\psi_1(x/|x|)}
\leq \lim_{|x|\to 0}\frac{\bar u_r(x)}{|x|^{a_{\mu_1}}\psi_1(x/|x|)}=\bar L_r.
\end{align*}}
Letting $r\to 0$, we complete the proof.
\end{pfn}

\section{Behavior of solutions to the semilinear problem}\label{sec:semi}

The $L^q$ and $L^{\infty}$ bounds of solutions to dipole-type linear
Schr\"odinger equations with properly summable potentials, derived in
Theorems \ref{t:bk} and \ref{t:bkn2}, allow us to obtain in  the semilinear
case analogous estimates.

\begin{Theorem}\label{t:bksemi}
  Let $\Omega$ be a bounded domain containing $0$, $a\in
  L^{\infty}({\mathbb S}^{N-1})$ such that $\Lambda_N(a)<1$, and
  $f:\Omega\times\R\to\R$ such that, for some positive constant $C$,
$$
\bigg|\frac{f(x,u)}u\bigg|\leq C\,\big(1+|u|^{2^*-2}\big)\quad\text{for a.e. }
(x,u)\in \Omega\times\R.
$$ 
Then, for any
  $\Omega' \Subset \Omega$ and 
for any weak
  $H^1(\Omega)$-solution  $u$ of \eqref{eq:72},
there holds  $\frac {u}{\varphi}\in
L^{\infty}(\Omega')$.
\end{Theorem}
\begin{pf}
Let $\Omega' \Subset \Omega$ and $u\in H^1(\Omega)$ be a
 weak 
  $H^1(\Omega)$-solution to \eqref{eq:72}. We set
$$
V(x)=\frac{f(x,u(x))}{\varphi^{2^*-2}(x)u(x)}
$$
and notice that $u\in L^{2^*}(\Omega)$ yields
$$
\int_{\Omega}\varphi^{2^*}(x)|V(x)|^{N/2}\,dx<+\infty.
$$
Hence Theorem \ref{t:bkn2} implies that $\frac {u}{\varphi}\in
L^{q}(\varphi^{2^*}\!\!,\Omega')$ for all $1\leq q<+\infty$. Since
$$
\int_{\Omega}\varphi^{2^*}(x)|V(x)|^{s}\,dx\leq {\rm
  const\,}\bigg(1+\int_{\Omega}\varphi^{2^*}(x)\bigg|
\frac{u(x)}{\varphi(x)}\bigg|^{(2^*-2)s}\,dx\bigg),
$$
we obtain that $V\in L^s(\varphi^{2^*}\!\!,\Omega)$ for all $s\geq
\frac{N-2}4$. The conclusion follows now from Theorem \ref{t:bk}.~\end{pf}

\noindent
\begin{pfn}{Theorem \ref{t:semi}} From Theorem \ref{t:bksemi}, it
   follows that $q(x)=\frac{f(x,u(x))}{u(x)}$ is such that $q\in
     L^{\infty}_{\rm loc}(\Omega\setminus\{0\})$ and
    $q(x)=O(|x|^{-(2-\e)})$ as $|x|\to0$ for some $\e>0$. Hence the
    conclusion follows from Theorem \ref{t:1}.
  \end{pfn}


\begin{thebibliography}{99}

\bibitem{AFP} B. Abdellaoui, V. Felli, I. Peral, 
{\it Existence and multiplicity for perturbations of an equation involving Hardy
inequality and critical Sobolev exponent in the whole ${\R}^N$,}
Adv. Differential Equations, 9 (2004), 481--508.

\bibitem{BrezisKato} H.~Br{\'e}zis, T.~Kato, {\it Remarks on the
    {S}chr\"odinger operator with singular complex potentials}, J.
  Math. Pures Appl. (9), 58 (1979), no.~2, 137--151.

\bibitem{CKN}
L.~Caffarelli, R.~Kohn,  L.~Nirenberg,
{\it First order interpolation inequalities with weights},
 Compositio Math.,  53 (1984), no.~3, 259--275.

\bibitem{CatrinaWang}
F.~Catrina, Z.-Q. Wang,
{\it On the {C}affarelli-{K}ohn-{N}irenberg inequalities: sharp
  constants, existence (and nonexistence), and symmetry of extremal
  functions},
 Comm. Pure Appl. Math.,  54 (2001), no.~2, 229--258.

\bibitem{DV2} V. De Cicco, M. A. Vivaldi, {\it A Liouville type theorem 
for weighted elliptic equations},  Adv. Math. Sci. Appl.,  9  (1999),  no. 1, 
183--207.

\bibitem{egnell} H. Egnell, {\it Elliptic boundary value problems with
    singular coefficients and critical nonlinearities,} Indiana Univ.
  Math., J. 38 (1989), no. 2, 235--251.


\bibitem{FMT} V. Felli, E.M. Marchini, S. Terracini, {\it On
    Schr\"odinger operators with multipolar inverse-square
    potentials,} preprint 2006, available online at
  \texttt{http://arxiv.org/abs/math.AP/0602209}.


\bibitem{FS3} V. Felli, M. Schneider, 
{\it A note on regularity of solutions to degenerate elliptic equations of
Caffarelli-Kohn-Nirenberg type,} 
Adv. Nonlinear Stud., 3 (2003), no. 4, 431--443.  

\bibitem{GP} J. Garc\'{\i}a Azorero, I. Peral, 
{\it Hardy Inequalities and some critical elliptic and parabolic problems,}
J. Diff. Equations, 144 (1998), no. 2, 441--476. 

\bibitem{Jan} E. Jannelli, 
{\it The role played by space dimension in elliptic critical problems,} 
J. Differential Equations, 156 (1999), no. 2, 407--426.  


\bibitem{leblond} J. M. L\'evy-Leblond, 
{\it Electron capture by polar molecules,} 
Phys. Rev., 153 (1967), no. 1, 1--4. 



\bibitem{murata} M. Murata, {\it Structure of positive solutions to
    $(-\Delta+V)u=0$ in $\R\sp n$,} Duke Math. J.  53 (1986), no. 4,
  869--943.

\bibitem{pinchover94} Y. Pinchover, {\it On positive Liouville
    theorems and asymptotic behavior of solutions of Fuchsian type
    elliptic operators,} Ann. Inst. H. Poincar\'e Anal. Non Linéaire, 11
  (1994), no. 3, 313--341.


\bibitem{reedsimon4} M. Reed, B. Simon, {\it Methods of modern
    mathematical physics. IV. Analysis of operators,} Academic Press,
  New York-London, 1978.

\bibitem{SafarovVassiliev} Yu. Safarov, D. Vassiliev, {\it The asymptotic
    distribution of eigenvalues of partial differential operators,}
    Translated from the Russian manuscript by the authors,
    Translations of Mathematical Monographs, 155. American
    Mathematical Society, Providence, RI, 1997.

\bibitem{simon2000} B. Simon, {\it Schr\"odinger operators in the twentieth century,}
  J. Math. Phys.,  41  (2000),  no. 6, 3523--3555.

\bibitem{SM} D. Smets, {\it Nonlinear Schr\"{o}dinger equations with
Hardy potential and critical nonlinearities,} 
Trans. AMS, 357 (2005), 2909--2938.


  \bibitem{terracini96} S. Terracini, {\it On positive entire solutions
      to a class of equations with singular coefficient and critical
      exponent,} Adv. Diff. Equa., 1 (1996), no. 2, 241--264.


  \end{thebibliography}
\end{document}